\documentstyle[12pt]{article} %\setlength{\textheight}{9in}
\newtheorem{prop}{Proposition}[subsection]
\newtheorem{dfn}[prop]{Definition}
\newtheorem{theo}[prop]{Theorem}

\newtheorem{rem}[prop]{Remark}
\newtheorem{coro}[prop]{Corollary}
\newtheorem{exam}[prop]{Example}

\newtheorem{lem}[prop]{Lemma}

\def\R{{\bf R }}

\def\Z{{\bf Z }}
\def\Q{{\bf Q }}

\def\L{ \Lambda }

\def\S{ \Sigma }

\def\ra{\rightarrow}

\begin{document}

\title{Rational Points of Bounded Height on Compactifications of
Anisotropic Tori}

\author{Victor V. Batyrev\thanks{Supported by DFG.} \\
\small Universit\"at-GHS-Essen, Fachbereich  6,  Mathematik \\
\small  Universit\"atsstr. 3,  45141  Essen, FRG  \\
\small e-mail: victor.batyrev@aixrs1.hrz.uni-essen.de \\
and \\
Yuri Tschinkel\thanks{Currently Junior Fellow of the
Harvard Society of Fellows.} \\
\small Harvard University, Department of Mathematics \\
\small 1 Oxford Street, Cambridge, MA 02138, USA  \\
\small e-mail: tschink@math.harvard.edu}

\date{}

\maketitle

\thispagestyle{empty}

\begin{abstract}
We investigate the analytic properties of  the zeta-function
associated with heights on equivariant compactifications
of anisotropic tori over number fields. This allows to verify
conjectures about the distribution of rational
points of bounded height.
\end{abstract}

\newpage

\tableofcontents

\newpage

{\Large {\bf Introduction}}

\bigskip

In this paper we prove new results on
the distribution of $K$-rational points of bounded height on
algebraic varieties $X$ defined over
a number field $K$ \cite{bat.man,franke-manin-tschinkel}.

Let ${\cal L}=(L, \| \cdot \|_v) $ be an ample metrized
invertible sheaf on $X$ with
a family $\{ \| \cdot \|_v \}$ of $v$-adic metrics and
$H_{\cal L}: X(K)\ra \R_{\ge
 0}$
the height function associated with $( L, \| \cdot \|_v )$.
We are interested in analytical properties of the
height
zeta-function defined by the series
\[ Z_{\cal L}(s) = \sum_{x \in X(K)} H_{\cal L}(x)^{-s}. \]
The investigation of $Z_{\cal L}(s)$
was  initiated by Arakelov and Faltings \cite{arakelov,faltings} who
considered the case when $X = S^d(C)$ is $d$-th symmetric power
of an algebraic curve $C$ of genus $g < d-1$.
It was shown in \cite{franke-manin-tschinkel} that zeta-functions
associated with heights on generalized flag varieties are Eisenstein series.

In this paper we consider the case when $X$ is an equivariant
compactification of an anisotropic torus $T$. It is easy
to show (cf. \ref{point.aniso}) that all $K$-rational points
of $X$ must be contained in $T(K)$. Therefore
\[ Z_{\cal L}(s) = \sum_{x \in T(K)} H_{\cal L}(x)^{-s}. \]
Our main idea for the computation of the zeta-function $Z_{\cal L}(s)$
is to use the {\em group structure} on $T(K)$
and to apply the Poisson formula in the following form:

{\em Let ${\cal G}$ be a locally compact topological abelian group
with a Haar measure $dx$, ${\cal H}
\subset {\cal G}$ a discrete  subgroup such that ${\cal G}/{\cal H}$
is compact,
$F \;: \; {\cal H} \rightarrow
{\bf R}$ a function on ${\cal H}$ which can be extended to an
$L^1$-function on ${\cal G}$. Let ${\rm  vol}(\cal G/\cal H)$ be
the $dx$-volume of the fundamental domain of ${\cal H}$ in
${\cal G}$. Denote by $\hat{F}$ the Fourier transform
of $F$ with respect to the Haar measure $dx$.
Suppose that  $\hat{F}\in L^1({\cal H}^{\perp})$.
Then
\[ \sum_{x \in \cal H} F(x) = \frac{1}{{\rm vol}(\cal G/\cal H) }
\sum_{\chi \in {\cal H}^{\perp}} \hat{F}(\chi) \]
where ${\cal H}^{\perp}$  is the group of  characters of
${\cal G}$ which are trivial on ${\cal H}$. }
\bigskip

We will apply this formula to the case when ${\cal G} = T({\bf A}_K)$
is the adele group of an anisotropic algebraic torus $T$,
${\cal H}$ is the
subgroup  $T(K)\subset T({\bf A}_K)$ of all $K$-rational points of $T$
and   $F(x,s) = H_{\cal L}^{-s}(x)$ for some ample
metrized line bundle ${\cal L}$ on the compactification of $T$.
The adelic method turns out to be very convenient, because
the height function $H_{\cal L}$ splits into the product of
local Weil functions
$H_{{\cal L},v} : T(K_v)\ra \R_{\ge 0}$. This allows to extend
$H_{\cal L}^{-s}$ to a function
on the whole adele group
${\cal G}$.
Moreover, the local Weil functions
$H_{{\cal L},v}$ can be chosen to be $T({\cal O}_v)$-invariant.
Therefore, the Fourier transform
\[ \int_{ T({\bf A}_K) } H_{\cal L}(x) \chi(x) d\mu
= \prod_{v \in {\rm Val}(K)} \int_{T(K_v)} H_{{\cal L},v}(x) \chi_v(x)
d\mu_v \]
equals zero unless the restriction
of $\chi_v$ on $T({\cal O}_v)$ is trivial for all $v
\in {\rm Val}(K)$.
The Fourier transform of $H_{\cal L}(x)^{-s}$ can be calculated
separately for each local factor $H_{{\cal L},v}(x)^{-s}$.
The analytic
properties of the Fourier transform of $H_{\cal L}(x)^{-s}$ can be investigated
by a method of Draxl \cite{drax1}.
\bigskip

In Section 1 we recall basic facts from theories of algebraic tori and
toric varieties ${\bf P}_{\Sigma}$ associated with
fans $\Sigma$ over arbitrary fields.
\medskip

In Sections 2,  following
ideas in \cite{franke-manin-tschinkel}, we define canonical families of
$v$-adic metrics on {\em all} $T$-linearized invertible sheaves
$L$ on ${\bf P}_{\Sigma}$ simultaneously. This allows us to construct
the complex height $H_{\Sigma}(x,\varphi)$
and the associated zeta-function $Z_{\Sigma}(\varphi)$ as
a function of  $\varphi$, where $\varphi$ represents an
element  of the complexified
Picard  group ${\rm Pic}({\bf P}_{\Sigma})_{\bf C}:
= {\rm Pic}({\bf P}_{\Sigma}) \otimes {\bf C}$. One obtains
the one-parameter zeta-function $Z_{\cal L}(s)$ via
restriction of $Z_{\Sigma}(\varphi)$
to the complex line $s \lbrack L \rbrack$, $s \in {\bf C}$.
\medskip

In Section 3 we recall some facts about characteristic functions of convex
cones.
\medskip

In Section 4 we investigate the analytic properties of zeta-functions in
order to obtain an asymptotic formula for the number of rational points
of bounded height. As one of our main results, we
prove the following refinement of the conjecture of Manin:
\medskip

{\em Let ${\bf P}_{\Sigma}$ be a smooth compactification of an anisotropic
torus $T$ $($notice that we do not need to assume
that ${\bf P}_{\Sigma}$ is a Fano variety$)$. Let $r$ be
the rank of ${\rm Pic}({\bf P}_{\Sigma,K})$.
Then there exist only a finite number $N( {\bf P}_{\Sigma},
{\cal K}^{-1}, B)$
of $K$-rational points  $x \in T(K)$ having the anticanonical
height $H_{{\cal K}^{-1}}(x) \leq B$.
Moreover,
\[ N({\bf P}_{\Sigma}, {\cal K}^{-1}, B) =
\frac{\Theta(\Sigma,K)}{(r-1)!}
\cdot B ( \log B)^{r-1}(1 + o(1)),\;\; B \ra \infty,  \]
where  the constant $\Theta(\Sigma,K)$ depends on:

$1.$ the  cone of effective divisors $\Lambda_{\rm eff}(\Sigma) \subset
{\rm Pic}({\bf P}_{\Sigma})_{\bf R}$;

$2.$ the Brauer group of  ${\bf P}_{\Sigma}$;

$3.$  the Tamagawa number $\tau_{\cal K}({\bf P}_{\Sigma})$
associated with
the metrized canonical sheaf on ${\bf P}_{\Sigma}$
as defined by E.Peyre in {\rm \cite{peyre}.} }
\medskip

We prove  the Batyrev-Manin conjecture \cite{bat.man} which
describes the asymptotic
for the number  of $K$-rational points $x \in T(K)$ such that
$H_{\cal L}(x) \leq B$ in terms
of $\lbrack {\cal L} \rbrack$ and
the geometry of the cone of effective divisors $\Lambda_{\rm eff}(\Sigma)$.
 Since for compactifications of anisotropic tori
  the cone of effective divisors
 $\Lambda_{\rm eff}(\Sigma) \subset
{\rm Pic}({\bf P}_{\Sigma})_{\bf R}$ is simplicial (cf. \ref{simp}),
our  situation is very close to the case of generalized flag
varieties considered in \cite{bat.man,franke-manin-tschinkel}.

We  observe that our results  provide first examples of  asymptotics
for the number of rational points of bounded height on
unirational varieties  of small dimension which are not rational
(nonrational anisotropic tori in dimension $3$
were constructed in \cite{kunyavskii}). Another new phenomenon
is the appearence of the Brauer group ${\rm Br} ({\bf P}_{\Sigma})$
in the asymptotic formulas.
\bigskip

{\bf Acknowledgements. } We are deeply grateful to Yu.I. Manin for
his permanent interest and strong encouragement of our work.

This paper was completed while the second author was visiting
the Arbeitsgruppe Zahlentheorie (MPI) in Berlin. He would like
to thank the Max-Planck-Institute for its hospitality and
ideal working conditions.

\section{Toric varieties over arbitrary fields}

\subsection{Algebraic tori}

Let $K$ be an arbitrary field, $\overline{K}$ the algebraic closure of
$K$, ${\bf G}_m(\overline{K}) = \overline{K}^*$ the multiplicative
group of $\overline{K}$.

Let $X$ be an arbitrary algebraic
variety over $\overline{K}$. Let $E/K$ be a finite extension such that $X$ is
defined over $E$. To stress this fact we sometimes will denote $X$ also
by $X_E$. The set of $E$-rational points of $X_E$ will be
denoted by $X_E(E)$.

\begin{dfn}
{\rm A linear algebraic group $T$
over  $K$ is
called  a {\em  $d$-dimen\-sio\-nal algebraic torus} if
its base extension $T_{\overline{K}} =
 T \times_{{\rm Spec}(K)} {\rm Spec }(\overline{K})$
is isomorphic to
$({\bf G}_m(\overline{K}))^d$.  }
\label{opr.tori}
\end{dfn}
We notice that an isomorphism
between ${T}$ and $({\bf G}_m(\overline{K}))^d$ is
always defined over a finite Galois extension $E$ of $K$.

\begin{dfn}
{\rm Let $T$ be an algebraic torus  over $K$. A
finite Galois extension
$E$ of $K$ such that $T_E = T \times_{{\rm Spec}(K)}
{\rm Spec}(E)$ is isomorphic to $({\bf G}_m(E))^d$
is called a {\em splitting field }  of  $T$.}
\end{dfn}

\begin{dfn}
{\rm We denote by $\hat{T} = {\rm Hom}\,( T, \overline{K}^*)$
the group of regular $\overline{K}$-rational
 characters of $T$.   For any subfield $E \subset
\overline{K}$
containing $K$, we denote by $\hat{T}_E$ the group of characters of $T$
defined over $E$. }
\end{dfn}
There is  well-known  correspondence between
Galois representations by integral matrices and algebraic tori
\cite{ono1,vosk}:

\begin{theo}
 Let $G = {\rm Gal }(E/K)$
be the Galois group of the splitting field $E$ of  a
$d$-dimen\-sional torus $T$. Then ${\hat{T}}$ is
a free abelian group of rank $d$ with a structure of $G$-module
defined by the natural representation
\medskip
\[ \rho \; :\; G \rightarrow {\rm Aut}(\hat{T}) \cong
{\rm GL}(d, {\bf Z}).\]
Every d-dimensional integral representation
of $G$ defines a $d$-dimen\-sional  algebraic torus over $K$ which splits
over $E$. One obtains a  one-to-one correspondence
between  $d$-dimensional
algebraic tori over
$K$ with the splitting field $E$ up to isomorphism, and $d$-dimensional
integral representations of $G$ up to equivalence.
\label{represent}
\end{theo}

\begin{rem}
{\rm The group $\hat{T}_K$ is a sublattice in
$\hat{T} \cong {\bf Z}^d$  consisting  of all $G$-invariant elements.}
\label{prop.char}
\end{rem}

\begin{dfn}
{\rm  An algebraic torus $T$ over $K$ is called
{\em anisotropic} if $\hat{T}_K$ has rank zero.}
\end{dfn}

\begin{exam}
{\rm
Let $f(z) \in K\lbrack z \rbrack$ be a separable
polynomial of degree $d$. Consider the $d$-dimensional $K$-algebra
\[ A(f) = K \lbrack z \rbrack / (f(z)). \]
Then the multiplicative group $A^*(f)$ is a $d$-dimensional algebraic
torus over $K$. This torus has the following properties:

(i) The rank of the group of characters of $A^*(f)$ is
equal to the number of irreducible components of ${\rm Spec}(A(f))$.

(ii) If  $f(z)$ splits in linear factors over some finite Galois extension $E$
of $K$, then
$A(f) \otimes_K E \cong E^n$, and $A^*(f)\otimes_K E \cong (E^*)^n$.
Thus, $E$ is a splitting field of $A^*(f)$.

(iii) Since the classes of $1, z, \ldots , z^{d-1}$ in $A(f)$ give
rise to a $K$-basis of the $d$-dimensional algebra $A(f)$, we can consider
$A^*(f)$ as a commutative subgroup in ${\rm GL}(d,K)$. Thus, the
determinant of the matrix defines  a regular $K$-character
\[{\cal N} \; : \; A^*(f) \rightarrow K^*. \]
We denote by $A^*_1(f)$ the $(d-1)$-dimensional algebraic torus which is
the kernel of ${\cal N}$.

(iv) The multiplicative group $K^*$ is a subgroup of $A^*(f)$ and the
restriction of  ${\cal N}$
to $K^*$ sends $x \in K^*$ to $x^d \in K^*$. The
factor-group $A^*(f) / K^*$ is a $(d-1)$-dimensional torus which
is isogeneous to $A^*_1(f)$. }
\label{exam.tori}
\end{exam}

\begin{exam}
{\rm Let $K'$ be a finite separable extension of $K$.
By primitive element theorem, $K' \cong A(f)$ for some
irreducible
polynomial $f(z) \in K \lbrack z \rbrack $.  Thus, we come to
a particular case of the previous example. In this case, ${\cal N}$ is
the norm $N_{K'/K}$, the
algebraic torus $A^*(f)$ is usually denoted by
\[ R_{K'/K} ({\bf G}_m) , \]
and the torus $A_1^*(f)$ is usually denoted by
\[ R_{K'/K}^1 ({\bf G}_m) . \]
Since ${\rm Spec}(K')$ is irreducible,
$R_{K'/K}^1 ({\bf G}_m)$ and $R_{K'/K} ({\bf G}_m)/ K^*$ are
examples of anisotropic tori.}
\label{exam.field}
\end{exam}

\subsection{Compactifications of split tori}

We recall standard facts about toric varieties over algebraically
closed fields \cite{danilov,demasur,oda}.

Let $M$ be a free abelian group of rank $d$ and $N = {\rm Hom}(M, {\bf Z})$ the
dual abelian group.

\begin{dfn}
{\rm A finite set $\Sigma$ consisting of  convex rational polyhedral
cones in $N_{\bf R} = N \otimes {\bf R}$ is called a {\em complete regular
$d$-dimensional fan} if the following conditions are satisfied:

(i) every cone $\sigma \in \Sigma$ contains $0 \in N_{\bf R}$;

(ii) every face $\sigma'$ of a cone $\sigma \in \Sigma$ belongs to $\Sigma$;

(iii) the intersection of any two cones in $\Sigma$ is a face of
both cones;

(iv) $N_{\bf R}$ is the union of cones from $\Sigma$;

(v)  every
cone $\sigma \in \Sigma$ is generated by a part of a ${\bf Z}$-basis of
$N$.\\
We denote by $\Sigma(i)$ the set of all $i$-dimensional cones in
$\Sigma$. For each cone $\sigma \in \Sigma$ we denote by
$N_{{\sigma}, \bf R}$ the minimal linear subspace containing $\sigma$. }
\label{def.fan}
\end{dfn}
\noindent
Every complete regular $d$-dimensional fan defines a smooth equivariant
compactification ${\bf P}_{\Sigma}$
of the split $d$-dimensional algebraic torus $T$. The variety
${\bf P}_{\Sigma}$ has the following two geometric properties:

\begin{prop}  The
toric variety ${\bf P}_{\Sigma}$ is the union of split  algebraic
tori $T_{\sigma}$ $($${\rm dim}\, T_{\sigma} = d - {\rm dim}\, \sigma$$)$:
\[ {\bf P}_{\Sigma} =  \bigcup_{ \sigma \in \Sigma } T_{\sigma}. \]
For each $k$-dimensional cone $\sigma \in \Sigma{(k)}$,
$T_{\sigma}$ is the kernel of a homomorphism $T \rightarrow
({\bf G}_m(\overline{K})^k)$ defined by a ${\bf Z}$-basis of
the sublattice $N \cap N_{{\sigma},{\bf R}} \subset N$.
\end{prop}
\noindent
Let $\check{\sigma}$ denote the cone in $M_{\bf R}$ which is
dual to $\sigma$.

\begin{prop} The toric variety ${\bf P}_{\Sigma}$ has a
$T$-invariant  open covering by affine subsets $U_{\sigma}$:
\[ {\bf P}_{\Sigma} = \bigcup_{ \sigma \in \Sigma} U_{\sigma} \]
where $U_{\sigma} = {\rm Spec}(\overline{K}) \lbrack M \cap \check{\sigma}
\rbrack$.
\end{prop}

\begin{dfn}
{\rm A continuous function $\varphi\; : \;
N_{\bf R} \rightarrow {\bf R}$ is called {\em $\Sigma$-piecewise linear}
if the restriction $\varphi_{\sigma}$ of $\varphi$ to every
cone $\sigma \in \Sigma$ is a linear function. It is
called {\em integral} if $\varphi(N) \subset {\bf Z}$.}
\end{dfn}

\begin{dfn}
{\rm For any integral $\Sigma$-piecewise linear function
$\varphi\; : \; N_{\bf R} \rightarrow {\bf R}$ and any cone $\sigma
\in \Sigma(d)$, we denote by
$m_{\sigma, \varphi}$ the restriction of $\varphi$ to $\sigma$
considered as an element in $M$.  We put $m_{\sigma', \varphi} =
m_{\sigma, \varphi}$ if $\sigma'$ is a face of a $d$-dimensional cone
$\sigma \in \Sigma$. }
\end{dfn}

\begin{dfn}
{\rm For any integral $\Sigma$-piecewise linear function
$\varphi\; : \; N_{\bf R} \rightarrow {\bf R}$, we define the
invertible sheaf $L(\varphi)$ as the subsheaf of the constant sheaf
of rational functions on ${\bf P}_{\Sigma}$ generated over
$U_{\sigma}$ by the element $- m_{\sigma,\varphi}$ considered as
a character of $T \subset {\bf P}_{\Sigma}$. }
\end{dfn}

\begin{rem}
{\rm The $T$-action on the sheaf of rational functions
restricts to  the subsheaf
$ L(\varphi)$ so that we can consider $ L(\varphi)$ as
a $T$-linearized line bundle over ${\bf P}_{\Sigma}$. }
\end{rem}

Denote by $e_1, \ldots, e_n$ the primitive integral generators
of all $1$-dimensional cones in $\Sigma$. Let $T_{i}$ $(i =1, \ldots, n)$
be the $(d-1)$-dimensional torus orbit corresponding to the cone
${\bf R}_{\geq 0}e_i \in \Sigma$ and  $D_i$ the Zariski closure of $T_i$ in
${\bf P}_{\Sigma}$. Define ${\bf D}({\Sigma}) \cong  {\bf Z}^n$
as the free abelian group  of  $T$-invariant Weil divisors on ${\bf
P}_{\Sigma}$ with the  basis $D_1, \ldots, D_n$.

\begin{prop}
The correspondence $\varphi \rightarrow L(\varphi)$
gives rise to an isomorphism between the group of
$T$-linearized line bundles on ${\bf P}_{\Sigma}$ and
the group $PL(\Sigma)$ of all
$\Sigma$-piecewise linear integral
functions on $N_{\bf R}$. There is the canonical isomorphism
\[ PL(\Sigma) \cong {\bf D}(\Sigma), \;\;
\varphi \mapsto (\varphi(e_1), \ldots, \varphi(e_n)). \]
The Picard group ${\rm Pic}({\bf P}_{\Sigma})$
is isomorphic to $PL(\Sigma)/M$ where elements of $M$
are considered as globally
linear integral functions on $N_{\bf R}$, so that we have the
exact sequence
\begin{equation}
 0 \rightarrow M \rightarrow {\bf D}(\Sigma) \rightarrow
{\rm Pic}({\bf P}_{\Sigma}) \rightarrow 0
\end{equation}
\end{prop}

\begin{dfn}
{\rm Let  $\Lambda_{\rm eff}(\Sigma)$ be the cone in
${\rm Pic}({\bf P}_{\S})$ generated by classes of effective divisors
on ${\bf P}_{\Sigma}$. Denote by ${\Lambda}_{\rm eff}^*(\Sigma)$ the
dual to $\Lambda_{\rm eff}(\Sigma)$ cone.}
\end{dfn}

\begin{prop}
$\Lambda_{\rm eff}(\Sigma)$ is generated by the classes
$\lbrack D_1 \rbrack, \ldots, \lbrack D_n \rbrack$.
\label{generators}
\end{prop}

{\em Proof.}  Any divisor $D$ on ${\bf P}_{\sigma}$ is
linearly equivalent to an integral
linear combination of $D_1, \ldots, D_n$. Assume that
$D = a_1 D_1 + \cdots + a_n D_n$ is effective. Then there exists a
rational function $f$ on ${\bf P}_{\Sigma}$ having no poles and zeros on
$T$ such that
\begin{equation}
 (f) + D \geq 0.
\label{effective}
\end{equation}
We can assume that $f$ is character of $T$ defined by an element $m_f \in M$.
Then the condition (\ref{effective}) is equivalent to
\begin{equation}
b_i =  \langle m_f, e_i \rangle + a_i \geq 0,\; i =1, \ldots, n
\end{equation}
Then $D' = b_1 D_1 + \cdots + b_n D_n$ is linearly equivalent to $D$.
So every effective class $\lbrack D \rbrack $ is a
non-negative integral linear combination of
$\lbrack D_1 \rbrack, \ldots, \lbrack D_n \rbrack$. \hfill $\Box$

\begin{prop}
Let $\varphi_{\Sigma}$ be the $\Sigma$-piecewise linear
integral function such that $\varphi(e_1) =
\cdots = \varphi(e_n)= 1$. Then $L(\varphi_{\Sigma})$ is
isomorphic to the $T$-linearized anticanonical line bundle
on ${\bf P}_{\Sigma}$.
\end{prop}

\begin{exam} {\sl Projective spaces}.
{\rm  Consider a $d$-dimensional fan $\Sigma$
whose $1$-dimensional cones are generated by $d+1$ elements
$e_1, \ldots , e_{d}, e_{d+1} = -(e_1 + \cdots + e_{d})$,
where  $\{ e_1, \ldots , e_{d} \}$
is a ${\bf Z}$-basis of ${d}$-dimensional lattice $N$, and
$k$-dimensional cones in $\Sigma$ are generated by all possible
 $k$-element subsets in $\{ e_1,\ldots, e_{d+1} \}$. Then the corresponding
compactification ${\bf P}_{\Sigma}$ of the  $d$-dimensional split torus is
${\bf P}^d$.}
\label{proj.space}
\end{exam}

\begin{rem}
{\rm It is easy to see that the combinatorial construction of toric
varieties ${\bf P}_{\Sigma}$
immediatelly extends to arbitrary fields $E$; i.e.,  using
a rational complete polyhedral fan $\Sigma$, one
can define the  toric variety ${\bf P}_{\Sigma,E}$
as the equivariant compactification of the split torus
$({\bf G}_m(E))^d$. }
\end{rem}

\subsection{Compactifications of nonsplit tori}

Let $T$ be a $d$-dimensional  algebraic torus over $K$
with a splitting field $E$ and
$G = {\rm Gal}\, (E/K)$. Denote by $M$  the lattice $\hat{T}$ and put
$N ={\rm Hom}\, (M, {\bf Z})$.  Let
$\rho^*$ be the integral representation of $G$ in ${\rm GL}(N)$
which is dual to $\rho$.
In order to construct  a projective compactification of $T$ over $K$,
we need a complete fan $\Sigma$ of  cones
having an additional combinatorial structure: an {\em
action of the Galois group } $G$  \cite{vosk}:

\begin{dfn}
{\rm  A complete  fan $\Sigma \subset N_{\bf R}$ is called
{\em $G$-invariant} if for any $g \in G$ and for any $\sigma \in \Sigma$, one
has  $\rho^*(g) (\sigma) \in \Sigma$. }
\label{opr.invar}
\end{dfn}

\begin{theo}
Let $\Sigma$ be a complete regular $G$-invariant fan in $N_{\bf R}$. Then
there exists a complete algebraic variety  ${\bf P}_{\Sigma,K}$ over $K$
such that its base extension ${\bf P}_{\Sigma,K} \otimes_{{\rm Spec}\, K}
{\rm Spec}\, E$ is isomorphic to
the toric variety ${\bf P}_{\Sigma,E}$ defined over $E$ by $\Sigma$.

Let $\Sigma^G$ be the subset of all $G$-invariant cones $\sigma \in \Sigma$.
Then
\[ {\bf P}_{\Sigma}(K) = \bigcup_{\sigma \in \Sigma^G} T_{\sigma}(K), \]
where $T_{\sigma}$ is the $(d - {\rm dim}\, \sigma)$-dimensional
algebraic torus over $K$  corresponding to the restriction of
the integral $G$-representation in ${\rm GL}(M)$ to
the sublattice  $(\hat{\sigma} \cap - \hat{\sigma}) \cap M \subset M$.
\label{decompos}
\end{theo}

Taking $G$-invariant elements in the short exact sequence
\[  0 \rightarrow M \rightarrow {\bf D}(\Sigma) \rightarrow
{\rm Pic}({\bf P}_{\Sigma,E}) \rightarrow 0 \]
we obtain the exact sequence
\begin{equation}
0 \rightarrow  M^G \rightarrow
{\bf D}(\Sigma)^G \rightarrow
{\rm Pic}({\bf P}_{\Sigma,E})^G \rightarrow H^1(G,M) \rightarrow 0
\label{short3}
\end{equation}

\begin{prop}
The group ${\rm Pic}({\bf P}_{\Sigma,E})^G$ is canonically
isomorphic to the Picard group ${\rm Pic}({\bf P}_{\Sigma,K})$. Moreover
$H^1(G, M)$  is the Picard group of $T$.
\end{prop}

\begin{coro}
The correspondence $\varphi \rightarrow L(\varphi)$
induces  an isomorphism between the group of
$T$-linearized invertible sheaves  on ${\bf P}_{\Sigma,K}$ and
the group $PL(\Sigma)^G$ of all
$\Sigma$-piecewise linear integral $G$-invariant
functions on $N_{\bf R}$. An invertible sheaf $ L$ on
${\bf P}_{\Sigma,K}$ admits a $T$-linearization if and only if
the restriction of $L$ on $T$ is trivial. In particular,
some tensor power of $ L$ always admits a $T$-linearization.
\end{coro}

\begin{coro}
Let $\Lambda_{\rm eff}(\Sigma,K)$ be the cone of effective divisors
of ${\bf P}_{\Sigma,K}$. Then $\Lambda_{\rm eff}(\Sigma,K)$ consists
of $G$-invariant elements in $\Lambda_{\rm eff}(\Sigma)$.
\end{coro}

\begin{coro}{\rm \cite{colliot1}}
Let ${\bf P}_{\Sigma,K}$ be a  compactification of an
anisotropic torus $T$.  Then
all $K$-rational points of ${\bf P}_{\Sigma,K}$ are contained in
$T$ itself.
\label{point.aniso}
\end{coro}

{\em Proof.} By \ref{decompos}, it is sufficient to prove
that for an anisotropic
torus $T$ defined by some Galois representation of $G$ in ${\rm GL}(M)$,
there is no  $G$-invariant cone $\sigma$ of
positive dimension in $\Sigma$.

Assume that a $k$-dimensional cone $\sigma$ with the generators $\{ e_{i_1},
\dots , e_{i_k} \}$ is $G$-invariant. Then $ e_{i_1} + \cdots + e_{i_k}$ is a
nonzero $G$-invariant integral vector in the interior of $\sigma$.
Hence the sublattice $N^G$ of $G$-invariant
elements in $N$ has  positive rank.
Thus $M^G \cong \hat{T}_K$ also has  positive rank. Contradiction.
\hfill $\Box$

\begin{prop}
Let ${\bf P}_{\Sigma,K}$ be a a compactification of an anisotropic torus
$T$. Then the cone of effective divisors
${\Lambda}_{\rm eff}({\Sigma,K})$ is
simplicial. The rank of the Picard group ${\rm Pic}({\bf P}_{\Sigma,K})$
equals to the number of $G$-orbits in $\Sigma(1)$.
\label{simp}
\end{prop}

{\em Proof.} Let $A_1({\bf P}_{\Sigma})$ be the group
of $1$-cycles on ${\bf P}_{\Sigma,E}$
modulo numerical equivalence. We identify $A_1({\bf P}_{\Sigma})$
with the dual to ${\rm Pic}({\bf P}_{\Sigma})$ group. Consider
the dual cone ${\Lambda}^*_{\rm eff}({\Sigma,K})$.
Since ${\Lambda}^*_{\rm eff}({\Sigma,K}) =
{\Lambda}^*_{\rm eff}({\Sigma})^G$, by
\ref{generators}, ${\Lambda}^*_{\rm eff}({\Sigma,K})$
consists of non-negative $G$-invariant
${\bf R}$-linear relations among primitive generators of $\Sigma(1)$.
Let
\[ \Sigma(1) = \Sigma_1(1) \cup \ldots \cup \Sigma_l(1)  \]
be the decomposition of $\Sigma(1)$ into a union of $G$-orbits.
Then every $G$-invariant linear relation among the primitive
generators $e_1 , \ldots , e_n$ of the
$1$-dimensional cones  has the  form
\[ \sum_{1 \leq i \leq l }
\lambda_i \left( \sum_{ \sigma_j \in \Sigma_i(1) }
e_j \right) = 0\;\;\;(\sigma_j = {\bf R}_{\geq 0} e_j). \]
For every $i$ ($1 \leq i \leq l$), the sum
\[   \sum_{ \sigma_j \in \Sigma_i(1) } e_j   \]
is a $G$-invariant element of the lattice $N$.
Since $T$ is anisotropic, $N^G = 0$ and all sums
$\sum_{ \sigma_j \in \Sigma_i(1) } e_j $ must be  equal to zero. These
integral relations  give rise to a ${\bf Z}$-basis
$r_1, \dots , r_l$ of the group of
integral linear relations among $e_1 , \ldots , e_n$.
 Thus
$A_1 ({\bf P}_{\Sigma})^G_{\bf R}$ is isomorphic to ${\bf Z}^l$ and the cone
$\Lambda_{\rm eff}^*({\Sigma,K})$
consists of nonnegative linear combinations of
$r_1, \dots , r_l$. So the cone $\Lambda_{\rm eff}({\Sigma,K})$ is also
an $l$-dimensional simplicial cone in
${\rm Pic}({\bf P}_{\Sigma,K}) \otimes {\bf R}$.
\hfill $\Box$

Below we consider several examples of compactifications of anisotropic tori.

\begin{exam}
{\rm  Consider a $d$-dimensional fan $\Sigma$ as in \ref{proj.space}.
It has a natural action
of the symmetric group $S_{d+1}$. Let $E$ be
a Galois extension of $K$ such that the Galois group
${\rm Gal}\, (E/K)$ is a subgroup of $S_{d+1}$ (for instance,
$E$ is a simple algebraic extension defined by an
$K$-irreducible polynomial $f$).
Then the action of $G$  on
$\Sigma$ defines  a ${d}$-dimensional toric variety
${\bf P}_{\Sigma,K}$
which over $E$ is isomorphic to $d$-dimensional
projective space; i.e. ${\bf P}_{\Sigma,K}$ is a
Severi-Brauer variety. In particular, if $E = K(f)$, then
${\bf P}_{\Sigma,K}$ is a compactification of the ${d}$-dimensional
anisotropic torus  $R_{E/K}({\bf G}_m)/K^*$.
Since ${\bf P}_{\Sigma,K}$ contains infintely many $K$-rational points,
${\bf P}_{\Sigma,K}$ is in fact  isomorphic  to  ${\bf P}^{d}$ over $K$. }
\label{exam.pn}
\end{exam}

\begin{exam}
{\rm A complete  fan $\Sigma$ is called {\em centrally  symmetric} if it is
invariant under the map $-Id$ of $N_{\bf R}$.

Let $\Sigma$ be a centrally symmetric $4$-dimensional
fan and let $E$ be an exension of $K$ of degree 2.
The  $d$-dimensional torus $T$
corresponding  to the integral representation of
${\rm Gal}(E/K) \cong {\bf Z}/2{\bf Z}$
by $Id$ and $-Id$ is isomorphic to the anisotropic torus
$(R^1_{E/K})^d$. The ${\bf Z}/2{\bf Z}$-invariant fan $\Sigma$
defines the compactification ${\bf P}_{\Sigma,K}$  of $(R^1_{E/K})^d$.  }
\label{two.aniso}
\end{exam}

\begin{exam}
{\rm Let $K'$ be a cubic extension of a number field $K$.
 We construct a smooth compactification
of the $2$-dimensional anisotropic $K$-torus $R^1_{K'/K}({\bf G}_m)$ as
follows.
Let $Y$ be the cubic surface in ${\bf P}^3$ defined by the equation
\[ N_{K'/K}(z_1,z_2,z_3) = z_0^3 \]
where $ N_{K'/K}(z_1,z_2,z_3)$ is the homogeneous cubic norm-form. Over
the algebraic closure $\overline{K}$ it is
isomorphic to the singular cubic surface $z_1 z_2 z_3 = z_0^3$.
The $3$ quadratic
singular points $p_1, p_2, p_3 \in Y_{\overline{K}}$ are
defined over a splitting field $E$ of  $R^1_{K'/K}({\bf G}_m)$ and the Galois
group $G = {\rm Gal}\, K'/K$ acts on  $\{ p_1, p_2, p_3 \}$ by permutations.
There exists a minimal simultaneous resolution $\psi\; : Y' \rightarrow Y$ of
singularities which is defined over $K$. By
contraction $\psi' : Y' \rightarrow X$ of  the proper pull-back of three
 $(-1)$ curves which are preimages of lines passing through the singular points
we obtain a
Del Pezzo surface $X$ of anticanonical degree $6$ which is  a smooth
compactification of the anisotropic torus $R^1_{K'/K}({\bf G}_m)$. }
\end{exam}

 Let $k$ be a finite
field of characteristic $p$ containing $q = p^n$
elements. Any finite extension $k'$ of $k$ is a cyclic Galois extension
and the group $G = {\rm Gal }(k'/k)$  is generated by the Frobenius
automorphism $\phi \; : \; z \rightarrow z^q$. By ~\ref{represent}, any
$d$-dimensional
algebraic torus $T$ over $k$ splitting over $k'$ is uniquely defined
by the  conjugacy class
in ${\rm GL} (d, {\bf Z})$
of the integral matrix
\[ \Phi = \rho ( \phi ) .\]

\noindent
The characteristic polynomial of the matrix $\Phi$ gives the following
formula obtained by T. Ono \cite{ono1} for the number of
$k$-rational points in $T$ :

\begin{theo}
Let $T$ be a $d$-dimensional  algebraic torus defined over
a finite field $k$. In the above notations, one has
the following formula for the number of $k$-rational points of $T$:
\[ {\rm Card} \lbrack T(k) \rbrack = (-1)^d{\rm det} ( \Phi - q\cdot Id ). \]
\label{fin.tori}
\end{theo}

\begin{prop}
Let ${\bf P}_{\Sigma}$ be a toric variety over a finite field $k$ defined by
a $\Phi$-invariant fan $\Sigma \subset N_{\bf R}$. For any $\Phi$-invariant
cone $\sigma \in \Sigma^G$, let  $M_{{\bf R}, \sigma} =
\check{\sigma} \cap (- \check{\sigma})$ be the maximal linear subspace
in the dual cone $\check{\sigma}
\subset M_{\bf R}$. Let $\Phi_{\sigma}$ be the restriction of
$\Phi$ on $M_{{\bf R},\sigma}$. Then
\[ {\rm Card}\lbrack {\bf P}_{\Sigma}(k)\rbrack = \sum_{\sigma \in \Sigma^G}
(-1)^{{\rm dim}\,\sigma} {\rm det}(\Phi_{\sigma} - q \cdot Id). \]
\label{point.var}
\end{prop}
\medskip

{\em Proof.} By \ref{decompos},
\[ {\bf P}_{\Sigma} (k) = \bigcup_{\sigma \in \Sigma^G} T_{\sigma}(k). \]
Observe that $k'$ is a splitting field for every algebraic torus
$T_{\sigma,k}$ defined by the $\rho$-action of $\Phi_{\sigma}$.
Now the statement follows  from \ref{fin.tori}.
\hfill $\Box$

\subsection{Algebraic tori  over  local and global
fields}

First we fix our notations.
Let ${\rm Val}(K)$ be the set of all valuations of a global field $K$.
For any $v \in {\rm Val}(K)$, we denote by $K_v$ the completion of
$K$ with respect to $v$.
Let $v$ be a non-archimedian absolute valuation  of a number  field $K$ and
$E$ a
finite Galois extension of $K$.
 Let ${\cal V}$ be an extension of $v$ to $E$,
 $E_{\cal V}$ the completion
of $E$ with respect to ${\cal V}$. Then
\[ {\rm Gal}(E_{\cal V}/ K_v ) \cong G_v \subset G, \]
where $G_{v}$ is the decomposition subgroup of
$G$ and $ K_v \otimes_K E \cong \prod_{{\cal V} \mid v} E_{\cal V}.  $
Let $T$ be an algebraic torus over $K$ with the splitting field $E$.
Denote by
$T_{K_v}=T\otimes K_v$.

\begin{dfn}
{\rm We denote the group of characters  $\hat{T}_{K_v} = M^{G_v}$ by
$M_v$ and the dual group ${\rm Hom}(\hat{T}_{K_v} , {\bf Z})
= N^{G_v}$  by $N_v$. }
\end{dfn}

\noindent
Let ${(K_v \otimes_K E)}^*$ and $E^*_{\cal V}$ be the multiplicative
groups of $K_v \otimes_K E$ and $E_{\cal V}$ respectively. One has
\[ T_{K_v} = {\rm Hom}_G (\hat{T}, {(K_v \otimes_K E )}^* ) =
{\rm Hom}_{G_v} (M, E_{\cal V}^* ). \]
Denote by ${\cal O}_{\cal V}$
the maximal compact subgroup in $E^*_{\cal V}$.
There is a short exact sequence
\[ 1 \rightarrow
{\cal O}_{\cal V} \rightarrow E_{\cal V}^* \rightarrow {\bf Z}
\rightarrow 1,  \;\; b \rightarrow {\rm ord} \mid b \mid_{\cal V}. \]
\noindent
Denote by $T({\cal O}_v)$ the maximal compact subgroup in $T(K_v)$.
Applying the functor
${\rm Hom}_{G_v}(M_v, * )$ to the short
exact sequence above, we obtain  the
short exact sequence
\[ 1 \rightarrow N_v \otimes {\cal O}_{\cal V} \rightarrow
N_v \otimes E_{\cal V}^* \rightarrow
N_v \rightarrow 1 \]
which induces an injective homomorphism
\[ \pi_v\; : \;   T(K_v) / T({\cal O}_v)
\hookrightarrow N_v = N^{G_v}. \]

\begin{prop} {\rm \cite{drax1}}
The homomorphism $\pi_v$ has finite cokernel. Moreover, $\pi_v$ is
an isomorphism if $E$ is unramified in $v$.
\label{pi-image}
\end{prop}

\begin{dfn}
{\rm Let $S$ be a finite subset of ${\rm Val}(K)$ containing
all archimedian and ramified  non-archimedian valuations of $K$.
We denote by $S_{\infty}$ the set of all archimedian
valuations of $K$ and put $S_0 = S \setminus S_{\infty}$. }
\end{dfn}

\medskip

Now we assume that $v$ is an archimedian absolute valuation, i.e., $K_v$ is
${\bf R}$ or ${\bf C}$. It is known that any torus over
${\bf R}$ is isomorphic to the product of some copies of  ${\bf C}^*$,
${\bf R}^*$, or $S^1 = \{ z \in {\bf C} \mid z \overline{z} =1 \}$. The
quotient $T(K_v) / T({\cal O}_v)$ is isomorphic to the ${\bf R}$-linear
space $N_v \otimes {\bf R}$.
The homomorphism $T(K_v) \rightarrow T(K_v) / T({\cal O}_v)$
is simply the logarithmic
mapping onto the Lie algebra of $T(K_v)$.  Hence, we obtain:

\begin{prop}
For any archimedian absolute valuation $v$, the
quotient $T(K_v) / T({\cal O}_v)$ can be canonically identified with
the real Lie
algebra of $T(K_v)$ embedded  in the $d$-dimensional ${\bf R}$-subspace
$N_{\bf R}$.
\end{prop}
\medskip

\begin{dfn}
{\rm Denote by $T({\bf A}_K)$ the adele group of $T$, i.e., the
restricted topological product
\[ {\prod_{v \in {\rm Val}(K)}} T(K_v) \]
consisting of all elements ${\bf t} = \{  t_v \} \in
\prod_{v \in {\rm Val}(K)}T(K_v)$ such that $t_v \in T({\cal O}_v)$ for
almost
all $v \in {\rm Val}(K)$.
Let
\[T^1({\bf  A}_K) = \{ {\bf t} \in T({\bf A}_K) \mid
\prod_{v \in {\rm Val}(K)}
\mid m(t_v) \mid_v = 1, \; {\rm for \; all}\; m \in \hat{T}_K \subset M  \}.
\]
We put also
\[ {\bf K}_T = \prod_{v \in {\rm Val}(K)} T({\cal O}_v), \]
 }
\end{dfn}

\begin{prop} {\rm \cite{ono1}}
The groups $T({\bf A}_K)$, $T^1({\bf  A}_K)$, $T(K)$,
${\bf K}_T$ have the following properties which are generalizations
of the corresponding properties of the adelization of ${\bf G}_m(K)$:

{\rm (i)} $T({\bf A}_K)/T^1({\bf  A}_K) \cong {\bf R}^k$, where $k$ is the rank
of $\hat{T}_K$;

{\rm (ii)} $T^1({\bf  A}_K)/T(K)$ is compact;

{\rm (iii)} $T^1({\bf  A}_K)/ {\bf K}_T \cdot T(K) $
is isomorphic to the direct product of a finite group
${\bf cl}(T_K)$ $($this
is an analog of the idele-classes group ${Cl}(K)$$)$ and
a connected compact abelian topological group which dimension
equals the rank $r'$ of the group of ${\cal O}_K$-units in $T(K)$
$($this rank equals $r_1 + r_2 -1$ for ${\bf G}_m$$)$;

{\rm (iv)} $W(T) = {\bf K}_T \cap T(K)$ is  a finite
group of all torsion elements in $T(K)$ $($this is the analog of the group
of roots of unity in ${\bf G}_m(K)$$)$.
\label{subgroups}
\end{prop}

The following theorem of A. Weil plays a fundamental role in the
definition  of adelic measures on algebraic varieties.

\begin{theo} {\rm \cite{peyre,weil1} } Let $X$ be an $d$-dimensional
 smooth algebraic variety over a global field $K$. Denote by
${\cal K}$ the canonical sheaf on $X$  with a family of local metrics
$\|\cdot \|_v $. Then these local metrics uniquely define
natural $v$-adic measures  $\omega_{{\cal K},v}$ on $X(K_v)$.

Let $U \subset X$ be a Zariski open
subset of $X$. Then for almost
all $v \in {\rm Val}(K)$ one has
\[ \int_{U({\cal O}_v)} \omega_{{\cal K},v} =
\frac{{\rm Card} \lbrack U(k_v) \rbrack}{q^d_v}, \]
where $k_v$ is the residue field of $K_v$ and $q_v = {\rm Card}
\lbrack k_v \rbrack$.
\end{theo}

\begin{rem}
{\rm We notice that the structure sheaf ${\cal O}_X$ of any
algebraic variety $X$ has a natural metrization defined by
$v$-adic valuations  of the field $K$.
If $X = {\cal G}$ is an algebraic group,
then there exists a natural
way to define a metrization of the canonical sheaf ${\cal K}$ on
${\cal G}$ by choosing a ${\cal G}$-invariant algebraic differential $d$-form
$\Omega$. Such a  form defines an isomorphism of ${\cal K}$ with
the structure sheaf ${\cal O}_{\cal G}$. We denote the corresponding
local measure on ${\cal G}(K_v)$ by  $\omega_{\Omega,v}$.}
\label{l.measures}
\end{rem}

Let $T$ be a $d$-dimensional torus over $K$ with a splitting field $E$.
Take a $T$-invariant
differential  $d$-form $\Omega$ on $T$ (it is unique up to a constant
from $K$). According to A. Weil (\ref{l.measures}),
we obtain a family of local measures $\omega_{\Omega,v}$ on $T(K_v)$.

\begin{dfn} {\rm \cite{ono1} Let
\[ L_S(s, T;E/K) = \prod_{v \not\in S} L_v(s, T ;E/K) \]
be the Artin $L$-function corresponding to the representation
\[ \rho \; :\; G= {\rm Gal}(E/K) \rightarrow {\rm GL}(M) \]
and a finite set $S \subset {\rm Val}(K)$ containing all
archimedian valuations and  all non-archimedian valuations of $K$ which are
ramified in $E$.
By definition,   $L_v(s,T;E/K) \equiv 1$ if $v \in S$. The numbers
\[ c_v = L_v(1, T; E/K) = \frac{1}{{\rm det}(Id - q^{-1}_v \Phi_v)},
\; v \not\in S \]
are called {\em canonical correcting factors} for measures
$\omega_{\Omega,v}$
($\Phi_v$ is the $\rho$-image of a local Frobenius element in
$G$).
 }
\end{dfn}
By \ref{fin.tori}, one has
\[ c_v^{-1}  = \int_{T({\cal O}_v)} \omega_{\Omega,v} =
\frac{{\rm Card} \lbrack T(k_v)\rbrack}{q^d_v}, \; \; v \not\in S. \]
Let $d\mu_v = c_v \omega_{\Omega,v}$. We put $c_v =1$ for $v \in S$. Since
\[ \int_{T({\cal O}_v)} d\mu_v = 1 \]
for  $v \not\in S$, the $\{ c_v \}$  defines
the canonical measure
\[ \omega_{\Omega,S} = \prod_{v \in {\rm Val}(K)} d\mu_v \]
on the adele group $T({\bf A}_K)$. By the product formula,
$\omega_{\Omega,S}$ does
not depend on the choice of $\Omega$.

Let $dx$ be the standard Lebesgue measure on $T({\bf A}_K)/T^1({\bf  A}_K)
\cong {\bf R}^k = M_{\bf R}^G$. There exists a unique Haar measure
$\omega^1_{\Omega,S}$ on $T^1({\bf A}_K)$ such that $\omega^1_{\Omega,S}dx =
\omega_{\Omega,S}$.

\begin{dfn}
{\rm The {\em Tamagawa number of } $T_K$ is defined as
\[ \tau(T_K) = \frac{b_S(T_K)}{l_S(T_K)} \]
where
\[ b_S(T_K) = \int_{T^1({\bf  A}_K)/T(K)} \omega^1_{\Omega,S} , \]
\[ l_S(T_K) = \lim_{s \rightarrow 1} (s-1)^k L_S(s, T; E/K). \]}
\label{tamagawa1}
\end{dfn}

\begin{rem}
{\rm Although the numbers $b_S(T_K)$ and $l_S(T_K)$ do depend on the choice
of the finite subset $S \subset {\rm Val}(K)$, the Tamagawa number
$\tau(T_K)$ does not depend on $S$.}
\end{rem}

\begin{theo} {\rm \cite{ono1,ono2} }
The Tamagawa number $\tau(T)$ of $T$ does not depend on the
choice of a splitting field $E$. It satisfies the following properties:

{\rm (i)} $\tau({\bf G}_m(K)) =1;$

{\rm (ii)} $\tau(T \times T') = \tau(T) \cdot \tau(T') $
where $T'$ and $T$ are tori over $K$;

{\rm (iii)} $\tau_K (R_{K'/K}(T)) = \tau_{K'}(T)$ for any
torus $T$ over $K'/K$.

Moreover, $\tau(T)$ is the ratio of two positive
integers
\[ h(T_K) = {\rm Card} \lbrack H^1(G, M) \rbrack \]
and $i(T_K) = {\rm Card}\lbrack {\rm III}(T) \rbrack$ where
\[ {\rm III}(T) = {\rm Ker}\, \lbrack
H^1(G, T(K)) \rightarrow \prod_{v} H^1(G_v, T(K_v)) \rbrack; \]
in particular, $\tau(T_K)$ is a rational number.
\label{tamagawa2}
\end{theo}

\begin{dfn}
{\rm
Let $\overline{T(K)}$ be the closure of $T(K)$ in $\prod_vT(K_v) $ in the
direct product topology. Define the {\em
obstruction group to weak approximation} as
\[ A(T)= \prod_v T(K_v)/\overline{T(K)}. \] }
\end{dfn}

\begin{theo} {\rm \cite{sansuc}}
Let ${\bf P}_{\S}$ be a complete smooth toric variety over $K$.
There is an exact sequence:
\[ 0\ra A(T) \ra Hom (H^1(G,{\rm Pic}({\bf P}_{\S,E})),\Q/\Z)\ra {\rm
III}(T)\ra
0.
\]
The group $H^1(G,{\rm Pic}({\bf P}_{\S,E}))$ is canonically isomorphic to
${\rm Br}({\bf P}_{\Sigma,K})/{\rm Br}(K)$, where
${\rm Br}({\bf P}_{\Sigma,K}) = H^2_{\rm et}({\bf P}_{\Sigma,K}, {\bf G}_m)$.
\label{weak}
\end{theo}

\begin{coro}
Denote by $\beta({\bf P}_{\Sigma})$ the cardinality of
$H^1(G,{\rm Pic}({\bf P}_{\S,E}))$. Then
\[ {\rm Card} \lbrack A(T) \rbrack =
\frac{\beta({\bf P}_{\Sigma})}{i(T_K)}. \]
\label{weak1}
\end{coro}

\section{Heights and their Fourier transforms}

\subsection{Complexified local Weil functions and heights}

A theory of heights on an algebraic variety $X$ defined over a number field $K$
 is
the unique functorial homomorphism from ${\rm Pic}(X) $ to equivalence
classes of functions
$X(K)\rightarrow {\bf R}_{\ge 0}$
which on metrized line bundles ${\cal L}$ is given by the formula
$$
H_{{\cal L}}(x) = \prod_v \| f(x)\|^{-1}_v
$$
where $f$ is a rational section of $L$ not vanishing in $x\in X(K)$.
Two functions are equivalent if
they differ by a bounded on $X(K)$ function.
For our purposes it will be convenient to extend these
notions to the complexified
Picard group ${\rm Pic}(X)\otimes {\bf C}$.

Let ${\bf P}_{\Sigma}$
be a compact toric  variety over a global field $K$.
We define a canonical
compact covering  of ${\bf P}_{\Sigma}(K_v)$
by compact subsets ${\bf C}_{\sigma,v} \subset U_{\sigma}(K_v)$.
For this purpose we identify lattice elements $m \in M$ with characters
of $T$ and define  the compact subset ${\bf C}_{\sigma,v} \subset
U_{\sigma}(K_v)$ as follows
\[ {\bf C}_{\sigma,v} = \{ x_v \in U_{\sigma}(K_v) \mid \| m(x_v) \|_v \leq 1\;
{\rm for}\; {\rm all}\; m  \in  M^{G_v}
\cap \check{\sigma} \}. \]

\begin{prop}
The compact subsets ${\bf C}_{\sigma,v}$ $(\sigma \in \Sigma)$ form a compact
covering of ${\bf P}_{\Sigma}(K_v)$ such that for any two cones
$\sigma, \sigma' \in \Sigma$ one has
\[ {\bf C}_{\sigma,v} \cap {\bf C}_{\sigma',v} = {\bf C}_{\sigma \cap
\sigma',v}. \]
\end{prop}

{\em Proof.} The last property of the compact subsets ${\bf C}_{\sigma,v}$
follows
immediatelly from their definition. Since the $T(K_v)$-orbit
of maximal dimension is dense in ${\bf P}_{\Sigma}(K)$, it is sufficient
to prove
that the compacts ${\bf C}_{\sigma,v}$ cover $T(K_v)$.

Let $x_v \in T(K_v)$. Denote by $\overline{x_v}$ the image of
 $x_v$
in $T(K_v) / T({\cal O}_v) \subset N_{\bf R}$.
By completeness of the fan $\Sigma$, the point $-\overline{x_v}$ is
contained in some cone $\sigma \in \Sigma$.
Hence  $x_v \in {\bf C}_{\sigma,v}$.
\hfill $\Box$
\bigskip

Now we define canonical metrizations of $T(K_v)$-linearized
line bundles on ${\bf P}_{\Sigma}(K_v)$.
\medskip

 Let $ L(\varphi)$ be a line bundle on ${\bf P}_{\Sigma}( K_v)$
corresponding to a $\Sigma$-piecewise linear integral $G_v$-invariant
function $\varphi$ on $N_{\bf R}$.

\begin{prop}
Let $f$ be a rational  section of $ L(\varphi)$. We define the $v$-norm
of $f$ at a point $x_v \in {\bf P}_{\Sigma}(K_v)$ as
\[ \| f(x_v) \|_v = \mid \frac{f(x_v)}{m_{\sigma,\varphi}(x_v)} \mid_v \]
where $\sigma$ is a cone in $\Sigma$ such that $x_v \in {\bf C}_{\sigma,v}$
and $m_{\sigma,\varphi} \in M$ is the restriction of $\varphi$ on $\sigma$.
Then this $v$-norm defines a $T({\cal O}_v)$-invariant
$v$-adic metric  on $ L(\varphi)$.
\end{prop}

{\em Proof.} The statement follows from the fact that
\[ \mid m_{\sigma,\varphi}(x_v) \mid_v = \mid m_{\sigma',\varphi}(x_v) \mid_v
\]
if $x_v \in {\bf C}_{\sigma,v} \cap {\bf C}_{\sigma',v}$.
\hfill $\Box$

A family of local metrics on all
$T$-linearized line bundles on ${\bf P}_{\Sigma}$ corresponding
to $\Sigma$-piecewise linear $G$-invariant
functions $\varphi \in PL(\Sigma)^G$   uniquely determines
a family of local Weil functions on $({\bf P}_{\Sigma})$
corresponding to $T$-invariant divisors
\[ D_{\varphi} = \varphi(e_1)D_1 + \cdots + \varphi(e_n)D_n. \]
We extend these local Weil functions to
the group of $T$-invariant Cartier divisors with {\em complex
coefficients} as follows.

\begin{dfn}
{\rm
A $T$-invariant ${\bf C}$-Cartier divisor is a formal linear combination
$D_s = s_1 D_1 + \cdots + s_n D_n$, with $s=(s_1,\ldots,s_n)\in {\bf C}^n$ or
equivalently
a complex piecewise linear function
$\varphi$ in $PL(\Sigma)^G_{\bf C}$ having the property
$\varphi(e_i) = s_i$ $(i = 1, \ldots ,n)$. }
\end{dfn}

\begin{dfn}
{\rm Let $\varphi \in PL(\Sigma)^G_{\bf C}$.
For any  point  $x_v \in T(K_v) \subset {\bf P}_{\Sigma}(K_v)$,
denote by $\overline{x_v}$ the image of
$x_v$ in $N_v$ (resp. $N_v\otimes {\bf R}$ for archimedian valuations),  where
$N_v$ is considered
as a canonical lattice  in the real space $N_{\bf R}$.
Define the {\em complexified
local Weil function}
$H_{\S,v}(x_v, \varphi)$ by the formula
\[H_{\S,v}(x_v,\varphi) = e^{\varphi(\overline{x_v})\log q_v }\]
where $q_v$ is  the cardinality of the residue field
$k_v$ of $K_v$ if $v$ is non-archimedian and $\log q_v = 1$ if
$v$ is archimedian. }
\end{dfn}

\begin{prop}
The complexified  local Weil function $H_{\S,v}(x_v, \varphi)$ satisfies the
following properties:

{\rm (i)} If $s_i = \varphi(e_i) \in {\bf Z}^n$ $( i =1, \ldots, n)$,
then  $H_{\S,v}(x_v,\varphi)$ is a classical local Weil
function $H_{L(\varphi),v}(x_v)$ corresponding to
a $T$-invariant Cartier divisor
\[ D_s = s_1 D_1 + \cdots + s_n D_n \]
on ${\bf P}_{\Sigma}$.

{\rm (ii)} $H_{\S,v}(x_v,\varphi)$ is  $T({\cal O}_v)$-invariant.

{\rm (iii)} $H_{\S,v}(x_v, \varphi + \varphi') =
 H_{\S,v}(x_v,\varphi)H_{\S,v}(x_v,\varphi')$.
\label{local.f}
\end{prop}

\begin{dfn}
{\rm Let $\varphi \in PL(\Sigma)^G_{\bf C}$. We define the {\em complexified
 height function on ${\bf P}_{\Sigma,K}$} by
\[ H_{\Sigma}(x,\varphi) = \prod_{v \in {\rm Val}(K)} H_{\S,v}(x,\varphi). \]}
\end{dfn}

\begin{rem}
{\rm Although all local factors $H_{\S,v}(x,\varphi)$ of
$H_{\Sigma}(x,\varphi)$ are functions on
$PL(\Sigma)_{\bf C}^G$, by the product formula,
the global complex height function $H_{\Sigma}(x,\varphi)$
depends only on the class of
$\varphi \in PL(\Sigma)_{\bf C}^G$ modulo
complex global linear $G$-invariant functions on $N_{\bf C}$,
i.e., $H_{\Sigma}(x,\varphi)$ depends only
on the class of $\varphi$ in ${\rm Pic}({\bf P}_{\Sigma,K})\otimes {\bf C}$. }
\end{rem}

\begin{dfn}
{\rm We define the zeta-function of the complex
height-function $H_{\Sigma}(x, \varphi)$ as
\[ Z_{\Sigma}(\varphi) = \sum_{x \in T(K)} H_{\Sigma}(x,-\varphi). \]}
\end{dfn}

\begin{rem}
{\rm One can see that the series $Z_{\Sigma}(\varphi)$
converges absolutely and uniformly in the domain
${\rm Re}(\varphi(e_j)) \gg 0$ for all $j$.
Since $H_{\Sigma}(x,\varphi)$ is the product of the local complex
Weil functions $H_{\S,v}(x,\varphi)$ and
$H_{\S,v}(x,\varphi) = 1$ for almost all $v$ ($x \in T(K)$),
we can immediately
extend $H_{\Sigma}(x,\varphi)$ to a function on the adelic group
$T({\bf A}_K)$.}
\end{rem}

\subsection{Fourier transforms of non-archimedian
heights}

\begin{dfn}
{\rm Let $\Sigma$ be a complete  regular fan of cones in $N$
whose 1-dimensional cones  are generated by
$e_1, \ldots, e_n$. We establish a
one-to-one correspondence between $e_1, \ldots, e_n$  and $n$
independent variables  $z_1, \ldots, z_n$.
The {\em Stanley-Reisner ring} $R(\Sigma)$ is defined as the factor of the
polynomial ring $A[z]= {\bf C} [z_1, \ldots, z_n ]$ by the ideal
$I(\Sigma)$ generated by all monomials
  $z_{i_1} \ldots z_{i_k}$ such that
   $e_{i_1}, \ldots, e_{i_k}$ are not generators of a $k$-dimensional
cone in $\Sigma$.}
\label{opr.stenley}
\end{dfn}

\begin{prop}
There is a natural identification  between the the elements of
the lattice $N$ and  the monomial {\bf C}-basis of the ring
$R(\Sigma)$.
\label{prop.basis}
\end{prop}

{\em Proof.}  Every integral point $x \in N$ belongs to the interior of
a unique cone $\sigma \in \Sigma$.
Let $e_{i_1}, \ldots, e_{i_k}$ be an integral basis of $\sigma$. Then
there exist positive integers $a_1, \ldots, a_k$ such that
\[ \overline{x} = a_1e_{i_1} + \cdots a_k e_{i_k}. \]
Therefore, $\overline{x}$ defines the
monomial $m(\overline{x}) = z_{i_1}^{a_1}  \cdots  z_{i_k}^{a_k}$.

By definition, $m(\overline{x}) \notin I(\Sigma)$.
It is clear that $I(\Sigma)$ has
a monomial {\bf C}-basis. Hence, we have
constructed a mapping  $\overline{x} \rightarrow m(\overline{x})$
from $N$ to
the monomial basis of $R(\Sigma)$.
It is easy to see that this mapping is  bijective.
\hfill $\Box$
\bigskip

Now choose a valuation  $v \not\in S$. Then we obtain
a cyclic subgroup  $G_v = \langle \Phi_v \rangle \subset G$
generated by a lattice automorphism
$\Phi_v \,:\, N \rightarrow  N$ representing
the local Frobenius element at place $v$.
Then $\Sigma(1)$ splits into  a disjoint union of $G_v$-orbits
\[  \Sigma(1) = \Sigma_1(1) \cup \cdots \cup \Sigma_l(1). \]
Let $d_j$ be the length of the $G_v$-orbit $\S_j(1)$.
One has
\[ \sum_{ i =1}^{l} d_j = n. \]

\begin{dfn}
{\rm Define the ${\bf Z}^l{\geq 0}$-grading of the
polynomial ring $A[z] = {\bf C} \lbrack z_1, \ldots, z_n \rbrack$
and the Stanley-Reisner ring $R(\Sigma)$
by the decomposition of the set of variables $\{z_1, \ldots, z_n\}$
into the disjoint union of $l$ sets $Z_1 \cup \cdots \cup Z_l$
which is induced by  the decomposition of
$\S(1)$ into $G_v$-orbits. The standard
${\bf Z}_{\geq 0}$-grading of the
polynomial ring $A[z] = {\bf C} \lbrack z_1, \ldots, z_n \rbrack$
and the Stanley-Reisner ring $R(\Sigma)$ will be
called the {\em total grading}.
}
\end{dfn}

\begin{dfn}
{\rm   We define the power series
$P(\Sigma, \Phi_v; t_1, \ldots, t_l)$ by the formula
\[ P(\Sigma,\Phi_v; t_1, \ldots, t_l) = \sum_{(i_1,\ldots, i_l) \in
{\bf Z}_{{\geq 0}}^l} ({\rm Tr}\,  \Phi_v^{i_1, \ldots, i_l})
t_1^{i_1} \cdots t_l^{i_l}, \]
where $\Phi_v^{i_1, \ldots, i_l}$ is the linear operator
induced by $\Phi_v$
on the
homogeneous $(i_1, \ldots, i_l)$-component of $R(\Sigma)$. }
\end{dfn}

\begin{prop}
One has
\[ P(\Sigma,\Phi_v; t_1, \ldots, t_l) = \frac{Q_{\S}(t_1^{d_1}, \ldots,
t_l^{d_l})}
{(1- t_1^{d_1})  \cdots (1 - t_l^{d_l}) } \]
where
$Q_{\S}(t_1^{d_1}, \ldots, t_l^{d_l})$ is a
polynomial in $t_1^{d_1}, \ldots, t_l^{d_l}$
having the total degree $n$ such its all
nonconstant monomials  have the total degree
at least $2$.
\label{p-function}
\end{prop}

{\em Proof.}
Since ${\rm dim} A \lbrack z \rbrack - {\rm dim} R(\Sigma) = n-d$,
there exists  the minimal
${\bf Z}^l_{\geq 0}$-graded free resolution
\[ 0 \ra F^{n-d} \ra \cdots \ra F^1 \ra F^0 = A \lbrack z \rbrack \ra
 R(\Sigma) \ra 0 \]
of the Stanley-Reisner ring $R(\Sigma)$ considered  as a module over the
polynomial ring $A \lbrack z \rbrack$.
Let $\overline{\Phi}_v^{\,i_1, \ldots, i_l}$ be the linear operator
on the homogeneous $(i_1, \ldots, i_l)$-component of $A\lbrack z
\rbrack$ induced by the action of $\Phi_v$ on $z_1, \ldots, z_n$. Then
\[ \sum_{(i_1,\ldots, i_l) \in
{\bf Z}_{{\geq 0}}^l} ({\rm Tr}\, \overline{\Phi}_v^{\,i_1, \ldots, i_l})
t_1^{i_1} \cdots t_r^{i_l} =  \frac{1}{\prod_{j = 1}^l (1 - t_j^{d_j})}. \]
We notice that ${\rm Tr}\, \overline{\Phi}_v^{\,i_1, \ldots, i_l}$ and
${\rm Tr}\, {\Phi}_v^{\,i_1, \ldots, i_l}$ can  be nonzero
only if the length $d_k$ of the $G_v$-orbit $\S_k(1)$
divides $i_k$ $( k =1, \ldots, l)$.
Therefore the polynomial $Q(t_1^{d_1}, \ldots, t_l^{d_l})$ is defined by
ranks and ${\bf Z}_{\geq 0}^l$-degrees of generators of
the free $A\lbrack z \rbrack$-modules $F^i$.
Notice that every monomial in $I(\Sigma)$ has the total degree at least $2$,
because every element $e_i \in \{ e_1, \ldots, e_n \}$ generates
a $1$-dimensional cone of $\Sigma$. So all
generators of $F^i$ $(i \geq 1)$ have the total degree at least $2$.
Therefore, the polynomial
$Q_{\S}$ has only monomials of the total degree at least $2$.
Since $R(\Sigma)$ is a Gorenstein ring,  $F^{n-d}$ is a free $A\lbrack z
\rbrack$-module of rank $1$ with a generator of the
degree $(d_1, \ldots, d_l)$. Therefore, $Q$ has the total degree
$n = d_1 + \cdots + d_l$. \hfill $\Box$
\bigskip

Let $\chi$ be a topological character of $T({\bf A}_K)$ such that
its $v$-component $\chi_v\, : \, T(K_v) \rightarrow S^1 \subset
{\bf C}^*$ is trivial on $T({\cal O}_v)$.
For each $ j \in \{ 1, \ldots, l\}$, we denote by
$n_j$ the sum of $d_j$ generators of all  $1$-dimensional
cones of the $G_v$-orbit $\S_j(1)$. Then $n_j$ is a $G_v$-invariant
element of $N$. By \ref{pi-image}, $n_j$ represents an element
of $T(K_v)$ modulo $T({\cal O}_v)$. Therefore, $\chi_v(n_j)$ is
well defined.

\begin{prop}
Let $v \not\in S$. Denote by $q_v$ the cardinality of the
finite residue field $k_v$ of $K_v$.  Then
for any local topological character
$\chi_v$ of $T(K_v)$,  one has
\[ \hat{H}_{\Sigma,v} (\chi, -\varphi) =
\int_{T(K_v)} H_{\S,v}(x_v,-\varphi) \chi_v(x_v)
 d\mu_v =  \]
 \[ \frac{Q_{\S}\left( \frac{\chi_{v}(n_1)}{{q_v}^{\varphi(n_1)}},
 \ldots, \frac{\chi_{v}(n_l)}{{q_v}^{\varphi(n_l)}} \right)}
{(1- \frac{\chi_{v}(n_1)}{{q_v}^{\varphi(n_1)}} )
\cdots (1 - \frac{\chi_{v}(n_l)}{{q_v}^{\varphi(n_l)}} ) }  \]
 if $\chi_v$ is trivial on $T({\cal O}_v)$, and
 \[ \int_{T(K_v)} H_{\S,v}(x_v,-\varphi) \chi_v(x_v)
 d\mu_v = 0 \]
 otherwise.
\label{integral.1}
\end{prop}

{\em Proof.} Since the local Haar measure $\mu_v$ is
$T({\cal O}_v)$-invariant, one has
\[ \int_{T(K_v)} H_{\S,v}(x_v, -\varphi) \chi_v(x_v) d \mu_v = \]
\[ = \sum_{\overline{x}_v
\in T(K_v)/T({\cal O}_v)} H_{\S,v}(\overline{x}_v, -\varphi)
\chi_v(\overline{x}_v) \int_{T({\cal O}_v)} \chi_v d\mu_v \]
where $\overline{x}_v$ denotes the image of $x_v$ in $T(K_v)/T({\cal O}_v)
= N_v$. Notice that $\int_{T({\cal O}_v)} \chi_v d\mu_v = 0$
if $\chi_v$ has nontrivial restriction on $T({\cal O}_v)$.
 By \ref{prop.basis},
there exists a natural identification between  $G_v$-invariant
elements of $N$ and $G_v$-invariant monomials in $R(\Sigma)$.
Since $\Phi_v$ acts by permutations on monomials in the homogeneous
$(i_1, \ldots, i_l)$-component of $R(\Sigma)$, the number
of $G_v$-invariant monomials in $R^{i_1, \ldots, i_l}(\Sigma)$ equals
${\rm Tr}\,\Phi_v^{i_1, \ldots, i_l}$.
Take a $G_v$-invariant element $\overline{x}_v \in N$ such that
$m(\overline{x}_v) \in R^{i_1, \ldots, i_r}(\Sigma)$.
Put  $i_k = d_k b_k$ $( k =1,
\ldots, l)$.
Then
\[ \varphi(\overline{x}_v) = b_1\varphi(n_1) + \cdots + b_l \varphi(n_l) \]
and
\[ \chi_v(\overline{x}_v) = \chi_{v}^{b_1}(n_1) \cdots \chi_{v}^{b_l}(n_l).
\]
This implies the claimed formula.
\hfill $\Box$
\bigskip

Let $A^*({\bf P}_{\Sigma}) = \bigoplus_{i =0}^d
A^i({\bf P}_{\Sigma})$ be the Chow ring
of ${\bf P}_{\Sigma,E_{\cal V}}$.
The groups $A^i({\bf P}_{\Sigma})$ have natural
$G_v$-action. Denote by $\Phi_v(i)$ the operator on $A^i({\bf P}_{\Sigma})$
induced by $\Phi_v$.

\begin{prop}
Denote by ${\bf 1}_v$ the trivial topological character of $T(K_v)$. Then
the restriction of
 \[ \int_{T(K_v)} H_{\S,v}(x_v,-\varphi) {\bf 1}_v(x_v) d \mu_v \]
to the line $s_1 = \cdots = s_r = s$ is equal to
 \[ L_v( s,T; E/K) \left( \sum_{k =0}^d
  \frac{{\rm Tr}\, \Phi_v(i)}{q_v^{ks}} \right). \]
  \label{loc-integ}
 \end{prop}

{\em Proof.} The Chow ring is the quotient of $R(\Sigma)$ by
a regular sequence \cite{danilov}. This gives the ${\bf Z}_{\geq 0}$-graded
(by the total degree) Koszul resolution having a $\Phi_v$-action:
\[ 0 \rightarrow \Lambda^d M \otimes R(\Sigma) \rightarrow
\cdots \rightarrow \Lambda^1 M \otimes R(\Sigma) \rightarrow
R(\Sigma) \rightarrow A^*({\bf P}_{\Sigma}) \rightarrow 0. \]
We apply the trace operator to the $k$-homogeneous component
of the Koszul complex, then we multiply the result by $1/q_v^{ks}$ and
take the sum over $k \geq 0$. By \ref{integral.1},
we have
\[ L_v^{-1}(s,T;E/K) \cdot
\int_{T(K_v)} H_{\S,v}(x_v,-\varphi) {\bf 1}_v(x_v) d \mu_v =
\sum_{k =0}^d \frac{{\rm Tr}\, \Phi_v(i)}{q_v^{ks}} \]
because
\[ \sum_{k =0}^d \frac{(-1)^k}{q_v^{ks}}
{\rm Tr}( \Lambda^k {\Phi_v}) =
{\rm det}(Id - q_v^{-s} \Phi_v) = L_v^{-1}(s,T;E/K). \]
\hfill $\Box$

\subsection{Fourier transforms of archimedian
heights}

\begin{prop}
Let $\chi_v(y) = e^{-2\pi i \langle x,y \rangle}$ be a topological character
of $T(K_v)$ which is trivial on $T({\cal O}_v)$. Then
the Fourier transform $\hat{H}_{\S,v}(\chi_v,-\varphi)$ of a
local archimedian Weil function
$H_{\S,v} (x,-\varphi)$ is a rational function in  $s_j = \varphi(e_j)$ for
${\rm Re}(s_j) > 0$.
\end{prop}

{\em Proof.} First we consider the case $K_v = {\bf C}$. Then
$T(K_v)/T({\cal O}_v) = N_{\bf R}$ and
\[ \frac{\hat{H}_{\S,v}(\chi_v(y),-\varphi)}{\int_{T({\cal O}_v)} d\mu_v}
 = \int_{N_{\bf R}} e^{-\varphi(x) -
2\pi i \langle x,y \rangle} dx =
\sum_{\sigma \in \Sigma(d)} \int_{\sigma} e^{-\varphi(x) -
2\pi i \langle x,y \rangle} dx \]
 where $dx=d\overline{x}_v$ is the standard measure
on $N_v\otimes \R\simeq \R^d $. On the other hand,
\[ \int_{\sigma} e^{-\varphi(x) -
2\pi i \langle x,y \rangle} dx = \frac{1}{\prod_{e_j \in \sigma}
(s_j + 2\pi  i \langle x,y \rangle)}. \]

If  $K_v = {\bf R}$, then the following simple
statements allows to repeat the arguments:

\begin{lem}
Let $\Sigma \subset N_{\bf R}$ be a complete regular $G$-invariant
fan of cones. Denote by $\Sigma^G \subset N_{\bf R}^G$
the fan consisting of $\overline{\sigma} = \sigma \cap N_{\bf R}^G$,
$\sigma \in \Sigma$. Then $\Sigma^G$ is again a complete
regular fan.
\end{lem}
 \hfill $\Box$
\medskip
\newline
The proof of the following proposition was suggested to us by W. Hoffmann.

\begin{prop}
Let ${\bf K} \subset {\bf C}^r$ be a compact such that ${\rm Re}( s_j) >
\delta$ for
all $(s_1, \ldots, s_r) \in {\bf K}$.  Then there exists a constant
$c({\bf K},\Sigma)$ such that
\[ \mid \hat{H}_{\S,v} (y,-\varphi) \mid \leq c({\bf K},\S) \sum_{\sigma \in
 \Sigma(d)}
\frac{1}{\prod_{e_k \in \sigma}
(1 + \mid \langle y, e_k \rangle \mid)^{1 + 1/d}}. \]

\label{l-estimation}
\end{prop}

{\em Proof.}
Let $f_1, \ldots, f_d$ be a basis of $M$. Put $x_i = \langle x, f_i \rangle$.
We denote by $y_1, \ldots, y_d$ the coordinates of $y$ in the basis
$f_1, \ldots, f_d$.
Let  $\varphi_i(x) = \frac{\partial}{\partial x_i}\varphi (x)$.
$\varphi_i(x)$ has a constant value $\varphi_{i,\sigma}$
in the interior of a cone $\sigma \in \Sigma(d)$.
\[ \hat{H}_{\S,v} (y,-\varphi) =  \int_{N_{\bf R}} e^{-\varphi(x)-2\pi
i<y,x>}dx =
 \frac{1}{2\pi iy_j} \int_{N_{\bf R}}
\frac{\partial}{\partial x_j} (e^{-\varphi(x)})
e^{-2\pi i<y,x>}dx \]
\[ = - \frac{1}{2\pi iy_j}
\int_{N_{\bf R}} \varphi_j(x) e^{-\varphi(x)-2\pi i<y,x>}dx \]
\[
= \frac{i}{2\pi y_j} \sum_{\sigma \in\Sigma(d)}
 \frac{\varphi_{j,\sigma}}{\prod_{e_k \in \sigma} (s_k+2\pi i<y,e_k>)} \]

Notice that $M_{\bf R}$ is covered by $d$ domains:
\[ V_j = \{ y = \sum_i y_i f_i \in M_{\bf R} \mid\;
\;  \mid  y_j \mid = \max_i \mid y_i \mid \}. \]
Let $\| y \|^2 = \sum_{i} y_i^2$. Then $\| y \| \leq \sqrt{d} \mid y_j \mid$
for $ y \in V_j$.  Then
\[ | \hat{H}_{\S,v} (y,-\varphi) |\leq
\frac{\sqrt{d}}{\|y\|}\sum_{\sigma \in \Sigma (d)}\frac{1}{
\prod_{e_k\in \sigma } |s_k+2\pi i<y,e_k>|} \]
for $y\in V_j$.
Furthermore, we obtain
\[ | \hat{H}_{\S,v} (y,-\varphi) |\leq
\frac{C'(\delta)}{1 + \|y\|} \sum_{\sigma \in \Sigma (d)}\frac{1}{
\prod_{e_k\in \sigma } (1 + \mid <y,e_k>\mid } \]
using the following obvious statement:

\begin{lem}
Assume that ${\rm Re}(s) > \delta > 0$. Then there exists
a positive constant $C(\delta)$ such that for all $t$ one has
$ C(\delta) (\mid s + 2\pi i t \mid) \geq 1 + \mid t \mid$.
\end{lem}

Since $\mid \langle y, e_k \rangle \mid \leq \| y \| \|e_k \|$, it follows
that there exist constants $c_{\sigma}$ such that
\[ c_{\sigma} (1+|<y,e_k>|)^d \geq \prod_{e_k\in \sigma }(1+|<y,e_k>|).
\]
Finally, we obtain
\[ \mid \hat{H}_{\S,v} (y,-\varphi) \mid \leq c_j(\delta, \Sigma)
\cdot \sum_{\sigma \in \Sigma(d)}
\frac{1}{\prod_{e_k \in \sigma}
(1 + \mid \langle y, e_k \rangle \mid)^{1 + 1/d}} \]
for all $y \in V_j$.
It remains to put $c({\bf K},\S) = \max_j c_j(\delta, \Sigma)$. \hfill $\Box$

\begin{coro}
Let $g \, : \, M_{\bf C} \ra {\bf C}$ be a continious function
such that
\[ \mid g(iy) \mid \leq  \|y\|^{\varepsilon}, \; \varepsilon < 1, \;\;
y \in M_{\bf R}. \]
Then
\[ \sum_{y \in O} g(iy) \hat{H}_{\S,v}(y, -\varphi) \]
is absolutely and uniformly convergent on ${\bf K}$ for any function $g(iy)$
for any lattice $O \subset M_{\bf R}$.
\label{lconver}
\end{coro}

\section{Characteristic functions of convex cones}

Let $V$ be an $r$-dimensional real vector space, $V_{\bf C}$ its
complex scalar extension, $\Lambda \subset V$
 a convex $r$-dimensional cone such that $\Lambda \cap - \L = 0
 \in V$. Denote by $\L^{\circ}$ the interior  of $\L$,
  ${\L}_{\bf C}^{\circ} =
{\L}^{\circ} + iV$
the complex  tube domain over ${\L}^{\circ}$, by
$V^*$  the
dual space,  by  ${\L}^* \subset V^*$  the dual to ${\L}$ cone and by
$dy$ a Haar measure on $V^*$.

\begin{dfn}
{\rm The {\em characteristic function of} ${\L}$ is defined as
the integral
\[  {\cal X}_{\L}(dy,u) =
\int_{{\L}^*} e^{- \langle u, y \rangle} dy, \]
where $u \in {\L}_{\bf C}$.  }
\end{dfn}

\begin{rem}
{\rm Characteristic functions of convex cones have been investigated
in the theory of homogeneous cones by M. K\"ocher, O.S. Rothaus, and
E.B. Vinberg \cite{koecher,vinberg,rothaus}.}
\end{rem}

\begin{rem}
{\rm We will be interested in characteristic functions of convex cones
${\L}$ in real spaces $V$ which have natural lattices $L \subset V$
of the maximal rank $r$. Let $L^*$ be the dual lattice in $V^*$, then
we can normalize the Haar measure $dy$  on $V^*$ so
that the volume of the fundamental domain
$V^*/L^*$ equals $1$. In this case the
corresponding characteristic function will be denoted simply by
${\cal X}_{\L}(u)$.  }
\end{rem}

\begin{prop} {\rm \cite{vinberg} }
Let $u \in {\L}^{\circ} \subset V$ be an interior point
of ${\Lambda}$. Denote by ${\L}^*_u(t)$ the convex $(r-1)$-dimensional
compact
\[ \{ y \in {\L}^* \mid \langle  u, y \rangle = t \} \]
We define the $(r-1)$-dimensional measure $dy_t'$ on ${\L}^*_u(t)$
in  such a way that for any function $f\; : \; V \rightarrow
{\bf R}$ with compact support one has
\[ \int_{{V}^*} f(y) dy = \int_{-\infty}^{+\infty}dt
\left(\int_{\langle u,y \rangle=t} f(y)dy_t'\right). \]
Then
\[ {\cal X}_{\L}(u) = (r-1)!\int_{{\L}^*_u(1)} dy_1'. \]
\end{prop}

The characteristic function  ${\cal X}_{\L}(u)$ has the following
properties \cite{rothaus,vinberg}:
\begin{prop}
{\rm (i)} If ${\cal A}$ is any invertible
linear operator on ${V}$, then
\[ {\cal X}_{\L} ({\cal A}u)  = \frac{{\cal X}_{\L}(u)}
{{\rm det}{\cal A}}; \]

{\rm (ii)} If ${\L}^{\circ} = {\bf R}^r_{\geq 0}$, $L = {\bf Z}^r \subset
{\bf R}^r$, then
\[ {\cal X}_{\L}(u) = (u_1 \cdots u_r)^{-1},  \;{\rm for }
 \;{\rm Re}(u_i) > 0 ; \]

{\rm (iii)} If $z \in {\L}^{\circ}$, then
\[ \lim_{z \rightarrow \partial {\L}} {\cal X}_{\L}(z) = \infty;  \]

{\rm (iv)} ${\cal X}_{\L}(u) \neq 0$ for all
$u \in {\L}_{\bf C}^{\circ}$.
\label{zeta.cone}
\end{prop}

\begin{prop}
If ${\L}$ is an $r$-dimensional finitely generated polyhedral cone,
then ${\cal X}_{\L}(u)$  is a rational function of degree $-r$.
In particular,  ${\cal X}_{\L}(u)$ has
a meromorphic extension to the whole complex space $V_{\bf C}$.
\end{prop}

{\em Proof.} It follows from Proposition \ref{zeta.cone}(i) that
${\cal X}_{\L}(\lambda u) =
{\lambda}^{-r} {\cal X}_{\L}(u)$. Hence ${\cal X}_{\L}(u)$ has
degree $-r$. In order to calculate ${\cal X}_{\L}(u)$,
we subdivide the dual cone ${\L}^*$ into
a union of simplicial subcones
\[ {\L}^* = \bigcup_{j} {\L}^*_j .\]
Then ${\L}$ is the intersection
\[ {\L} = \bigcap_j {\L}_j. \]
For ${\rm Re}(u) \in \bigcap_j {\L}_{j}^{\circ}$, one has
\[ {\cal X}_{\L}(u) = \sum_j {\cal X}_{{\L}_j}(u). \]
By Proposition \ref{zeta.cone}(i),(ii), every
function ${\cal X}_{{\L}_j}(u)$ is rational. \hfill $\Box$

\begin{dfn}
{\rm Let $X$ be a smooth proper algebraic variety. Denote by
$\L_{\rm eff} \subset {\rm Pic}(X)_{\bf R}$  the cone
generated by classes of effective divisors on $X$.
Assume that the anticanonical class $ \lbrack  {\cal K}^{-1}
\rbrack \in  {\rm Pic}(X)_{\bf R}$
is contained in the interior of $\L_{\rm eff}$. We define
the constant $\alpha(X)$ by
\[ \alpha(X) = {\cal X}_{\L_{\rm eff}}( \lbrack {\cal K}^{-1}
\rbrack). \]
}
\end{dfn}

\begin{coro}
If ${\L}_{\rm eff}$ is a finitely generated polyhedral cone,
then $\alpha(X)$ is a rational number.
\end{coro}

\begin{exam}
{\rm Let ${\bf P}_{\Sigma,K}$ be a smooth compactification of
an anisotropic torus $T_K$.
By \ref{simp}, $\L_{\rm eff} \subset
{\rm Pic}({\bf P}_{\Sigma,K}) \otimes {\bf R}$ is a
simplicial cone. Using \ref{zeta.cone} and the exact sequence
\[ 0 \ra PL(\Sigma)^G \ra {\rm Pic}({\bf P}_{\Sigma,K})
\ra H^1(G,M) \ra 0 \]
we obtain
\[ {\cal X}_{\L_{\rm eff}}(u) = \frac{1}{h(T_K) u_1
\cdots u_r}, \]
where $u = \varphi$, $\varphi(e_j) = u_j$ $(j =1, \ldots l)$.
In particular,
\[ \alpha({\bf P}_{\S}) = \frac{1}{h(T_K)}. \]
}
\end{exam}

\begin{exam}
{\rm Consider an example of a non-simplicial  cone of Mori $\L_{\rm eff}$
in $V = {\rm Pic}(X)_{\bf R}$ where $X$ is a Del Pezzo
surface  of anticanonical degree 6.
The cone ${\L}$ has 6 generators corresponding to exceptional curves
of the first kind on $X$. We can construct  $X$ as the blow up of
3 points $p_1, p_2, p_3$ in general position on ${\bf P}^2$.
The  exceptional curves are $C_1, C_2, C_3,
C_{12}, C_{13}, C_{23}$, where $C_{ij}$ is the proper pullback of the line
joining $p_i$ and $p_j$.

If $u = u_1 [C_1] + u_2 [C_2] + u_3 [C_3] + u_{12}[C_{12}] +
u_{13}[C_{13}] + u_{23} [C_{23}] \in \L_{\rm eff}^{\circ}$, then

\[ {\cal X}_{\L_{\rm eff}}(u) =
\frac{ u_1 + u_2 + u_3 + u_{12} + u_{13} + u_{23} }
{(u_1 + u_{23}) (u_2 + u_{13})(u_3 + u_{12})(u_1 + u_2 + u_3 )
(u_{12} + u_{13} + u_{23})} \]
and
\[ \alpha(X) = 1/12. \]}
\end{exam}

\begin{prop}
Assume that ${\L}$ is a finitely generated
polyhedral cone and ${\L} \cap -{\L} = 0$.
Let $p_0$ and $p_1$ be two points in ${E}$ such that $p_1 \not\in {\L}$
and $p_0 \in {\L}^{\circ}$. Let $t_0$ be a positive real number such
that $t_0p_0 + p_1 \in \partial {\L}$. We define a meromorphic function
in one complex variable $t$ as
\[ Z(p_0,p_1, t) = {\cal X}_{\L}(tp_0 + p_1 ). \]
Let $k$ be
the codimension of the minimal face of ${\L}$ containing
$t_0p_0 + p_1$. Then
the rational function $Z(p_0,p_1, t)$ is  analytic
for ${\rm Re}(t) > t_0$ and
it  has a pole of order $k$ at $t = t_0$.
\end{prop}

{\em Proof.} As in the proof of the previous statement, we can subdivide the
dual cone ${\L}^*$ into simplicial subcones ${\L}^*_j$ such
that $t_0p_0 + p_1 \in \partial {\L}_1$ and
$t_0p_0 + p_1 \not\in \partial {\L}_j$
$(j >1)$. It suffices now to apply Proposition \ref{zeta.cone}(i),(ii) to
${\L}_1$. \hfill $\Box$

\begin{coro}
Assume that ${\L}$ is only locally polyhedral at the
point $t_0p_0 + p_1$ and $k$ is
the codimension of a minimal polyhedral face of ${\L}$ containing
$t_0p_0 + p_1$. Then

{\rm (i)} $Z(p_0,p_1, t)$ is an analytical function for ${\rm Re}(t) > t_0$.

{\rm (ii)} $Z(p_0,p_1, t)$ has  meromorphic continuation to some
neigbourhood of $t_0$.

{\rm (iii)} $Z(p_0,p_1, t)$ has a pole of order $k$ at $t = t_0$.
\end{coro}

\section{Distribution of rational points}

\subsection{The method of Draxl}

Let $\Sigma$ be a $G$-invariant regular fan,
$\Sigma(1) = \Sigma_1(1) \cup \cdots \cup \Sigma_r(1)$ be
the decomposition of $\Sigma(1)$ into $G$-orbits. We choose a
representative $\sigma_j$  in each $\Sigma_j(1)$ $( j =1, \ldots, r)$.
Let $e_j$ be the primitive integral generator of $\sigma_j$,
$G_j \subset G$ be the stabilizer of $e_j$. Denote by $k_j$ the
length of $G$-orbit of $e_j$, and by $K_j \subset E$ the subfield of
$G_j$-fixed elements. Then $k_j = \lbrack K_j : K \rbrack$
$(j =1, \ldots, r)$.

Consider the $n$-dimensional torus
\[ T' := \prod_{j =1}^r R_{K_j/K}({\bf G}_m). \]
Notice that the group ${\bf D}(\Sigma)$ can be identified with the $G$-module
$\hat{T}'_K$.  The homomorphism of $G$-modules
$ M \rightarrow {\bf D}(\Sigma)$
induces the  homomorphism $T' \ra T$
and a map
$$
\gamma:\, \prod_{j=1}^r
{\bf G}_{m}({\bf A}_{K_j})/ {\bf G}_{m}(K_j) \ra T({\bf A}_{K})/T(K)
$$
We get a map of characters
$$
\gamma^*:\,  (T({\bf A}_{K})/T(K))^*\ra \prod_{j =1}^r
({\bf G}_{m}({\bf A}_{K_j})/{\rm G}_m(K_j))^*.
$$

\begin{rem}{\rm
The kernel of $\gamma^*$ is dual to the obstruction group to weak approximation
$A(T) $ defined above. }
\label{obstr}
\end{rem}

Let
\[ \chi \;: \; T({\bf A}_K) \rightarrow S^1 \subset  {\bf C}^* \]
be a topological character which is trivial on
$T(K)$.
Then $\chi \circ \gamma$ defines Hecke characters of the
idele groups
\[ \chi_j \; :\;
{\bf G}_m({\bf A}_{K_j}) \rightarrow S^1 \subset  {\bf C}^*. \]
If $\chi$ is trivial on ${\bf K}_T$, then all characters $\chi_j$
$(j =1, \ldots, r)$ are  trivial on the maximal
compact subgroups in ${\bf G}_m({\bf A}_{K_j})$.
We denote by  $L_{K_j}(s,\chi_j)$ the Hecke $L$-function corresponding to the
character $\chi_j$.
The following statement is well-known:

\begin{theo}
The function $L_{K_j}(s,\chi_j)$ is holomorphic in the whole plane
 unless $\chi_j$ is trivial. In the later case,
$L_{K_j}(s,\chi_j)$ is holomorphic for
${\rm Re}(s) >1$ and has a meromorphic extension to the complex plane  with
a pole of order $1$  at $s = 1$.
\end{theo}

We come to the main  statement which describes the analytical properties
of the Fourier transform of height functions.

\begin{theo}
Define  affine complex coordinates
$\{s_1, \ldots, s_r \}$ on the vector space
$PL(\Sigma)_{\bf C}^G$ by  $s_j = \varphi(e_j)$ $( j =1, \ldots, r)$.
Then the Fourier transform
$\hat{H}_{\S}(\chi, -\varphi)$ of the complex height
function $H_{\S}(x,-\varphi)$  is always an analytic function for
${\rm Re}(s_j) > 1$ $(1 \leq j \leq l)$, and
\[ \hat{H}_{\S}(\chi,-\varphi)
\prod_{i=1}^r L^{-1}_{K_j}(s_j,\chi_j) \]
has an analytic extension to the domain ${\rm Re}(s_j) > 1/2$
$(1 \leq j \leq r)$.
\label{dmethod}
\end{theo}

{\em Proof.} The idea of the proof is essentialy due to Draxl
\cite{drax1}. We have  the Euler product
\[ \hat{H}_{\S}(\chi,-\varphi) =
\prod_{v \in {\rm Val}(K)} \hat{H}_{\S,v}(\chi_v,-\varphi) \]
In order to prove the above properties of
$\hat{H}_{\S}(\chi, -\varphi)$, it is
sufficient to investigate the product
\[ \hat{H}_{\S,S}(\chi,-\varphi) = \prod_{v \not\in S}
\hat{H}_{\S,v}(\chi_v,-\varphi). \]

Choose a valuation  $v \not\in S$. Then we obtain
a cyclic subgroup  $G_v = \langle \Phi_v \rangle \subset G$
generated by a lattice automorphism
$\Phi_v \,:\, N \rightarrow  N$ representing
the local Frobenius element at place $v$.
Let $l$ be the number of $G_v$-orbits in $\S(1)$.
By \ref{integral.1},
\[ \hat{H}_{\S,S}(\chi, -\varphi) = \prod_{v \not\in S}
P\left( \Sigma, \Phi_v; \frac{\chi_{v}(n_1)}{{q_v}^{\varphi(n_1)}},
 \ldots, \frac{\chi_{v}(n_l)}{{q_v}^{\varphi(n_l)}} \right). \]
By \ref{p-function}, we have
\[ \hat{H}_{\S,v}(\chi, -\varphi) =
 \frac{Q_{\S}\left( \frac{\chi_{v}(n_1)}{{q_v}^{\varphi(n_1)}},
 \ldots, \frac{\chi_{v}(n_l)} {{q_v}^{\varphi(n_l)}} \right)}
 {(1- \frac{\chi_{v}(n_1)}{{q_v}^{\varphi(n_1)}} )
\cdots (1 - \frac{\chi_{v}(n_l)}{{q_v}^{\varphi(n_l)}} )} .   \]
Moreover,
\[ \prod_{v \not\in S}
Q_{\S}\left( \frac{\chi_{v}(n_1)}{{q_v}^{\varphi(n_1)}},
 \ldots, \frac{\chi_{v}(n_l)}{{q_v}^{\varphi(n_r)}} \right) \]
is an absolutely convergent Euler product for
${\rm Re}(s_j) > 1/2$ $( j =1, \ldots, l)$.

It remains to show the relation between
\[  \left( 1 - \frac{\chi_{v}(n_1)}{{q_v}^{\varphi(n_1)}} \right)^{-1} \cdots
\left( 1 - \frac{\chi_{v}(n_l)}{{q_v}^{\varphi(n_l)}} \right)^{-1} \]
and local factors of the product of the Hecke $L$-functions
\[ \prod_{j =1}^r L_{K_j}(s_j, \chi_j).  \]

For this purpose, we compare two decompositions of $\S(1)$ into
the disjoint union of $G_v$-orbits  and $G$-orbits.
Notice that for every $j \in \{ 1, \ldots, r \}$, the $G$-orbit
$\Sigma^j(1)$ decomposes into  a disjoint union of $G_v$-orbits
\[  \Sigma_j(1) = \Sigma_{j1}(1) \cup \cdots \cup \Sigma_{j l_j}(1). \]
Let $d_{ji}$ be the length of the $G_v$-orbit $\S^j_i(1)$; i.e., we put
$\{ d_{ji} \} = \{ d_1, \ldots, d_l \}$.
One has
\[ \sum_{ i =1}^{l_j} d_{ji} = k_j. \]
and
\[ l = \sum_{i =1}^r l_j. \]

On the other hand, $l_j$ is the number of different
valuations ${\cal V}_{j1}, \ldots, {\cal V}_{jl_j} \in {\rm Val}(K_j)$
over  of $v \in {\rm Val}(K)$. Let $k_v$ be the residue
field of $v \in {\rm Val}(K)$, $k_{{\cal V}_{ji}}$ the residue field
of ${\cal V}_{ji} \in {\rm Val}(K_j)$. Then
\[ d_{ji} =  \lbrack k_{{\cal V}_{ji}} : k_v \rbrack. \]

We put also $\{ n_1, \ldots, n_l \} =
\{ n_{ji} \}$, where $n_{ji}$ denotes the sum of $d_{ji}$ generators
of all $1$-dimensional cones of the $G_v$-orbit $\S_{ji}(1)$.
Therefore, $\chi_v(n_{ji})$ is the
${\cal V}_{ji}$-adic component
of the Hecke character $\chi_j$. Hence
\[ \prod_{i =1}^{l_j} \left( 1 - \frac{
\chi_{v}(n_{ji})}{{q_v}^{\varphi(n_{ji})}} \right)^{-1} \]
equals the product of the local factors
\[ \prod_{{\cal V}_{ji}}  \left( 1 -
\frac{\chi_{{\cal V}_{ji}}}{q_{{\cal V}_{ji}}^{\varphi(n_{ji})}}
\right)^{-1} \]
of the Hecke $L$-function $L_{K_j}(s_j, \chi_j)$.
\hfill $\Box$

\subsection{The meromorphic extension of
$Z_{\Sigma}(\varphi)$}

\begin{theo}
For any $\varepsilon > 0$ there exists a $\delta >0$ and a
constant $c({\varepsilon})$ such
that
\[ \mid L_K(s,\chi) \mid \leq c(\varepsilon)(\mbox{\rm Im}(s))^{\varepsilon}\;
\; \mbox{\rm for}\; u= \mbox{\rm Re}(s) > 1 - \delta  \]
for every Hecke $L$-function $L_K(s,\chi)$  with a
nontrivial nonramified character $\chi$.
\label{estim}
\end{theo}

{\em Proof.}  We use the following standard statement
based on the Phragm\'en-Lindel\"of principle:

\begin{lem} {\rm ( \cite{titchmarsh}, p.181)}
Let $f(s) $ be a single valued analytic function in the strip
$u_1 \le {\rm Re}(s)\le u_2$ satisfying the  conditions:

{\rm (i)} $|f(u + it )| < A_0 \exp ( e^{C|t|})$ for some
real constants $A_0 >0$ and $ 0 < C < \pi/(u_2 - u_1)$;

{\rm (ii)} $| f(u_1 + it) | \leq A_1 |t|^{a_1}$,
 $| f(u_2 + it) | \leq A_2 |t|^{a_2}$ for some
 constants $a_1$, $a_2$.
Then for all $u_1 \leq u \leq u_2$, we have the estimate
\[ | f(u + it) | \leq A_3 |t|^{a(u)} \]
where
\[ a(u) =a_1  \frac{u_2-u}{u_2-u_1} + a_2 {\frac{u-u_1}{u_2-u_1}}.\]
\label{FL}
\end{lem}

Choose a sufficiently small $\delta_1$. Then $| L_K(s,\chi) |$ is
bounded by $A_2(\delta_1)=
 \zeta_K(1+\delta_1)$ for ${\rm Re}(s) = 1 + \delta_1$.
Consider  the functional equation  $L_K(s,\chi) =
C(s)L_K(1-s,\overline{\chi})$.
Since $\chi$ (as well as  $\overline{\chi}$) is unramified, the function
$C(s)$ depends only on the field $K$.
Using standard estimates for $\Gamma$-factors in  $C(s)$, we obtain
$| L_K(s, \chi) | <  A_1 ({\rm Im}(s))^{a_1}$ for ${\rm Re}(s) =
- \delta_1$ and some sufficiently large explicit constants $A_1(\delta_1),
a_1(\delta_1)$.

We apply Lemma \ref{FL} to the $L$-function
$L_K(s,\chi)$ where $u_1 = - \delta_1$,
$u_2 = 1 + \delta_1$, and $a_2 = 0$.
Then for $1 - \delta < {\rm Re}(s) < 1 + \delta_1$, one has
\[ | L_K(s,\chi) | \leq A_3(\delta,\delta_1)({\rm Im}
(s))^{a_1\frac{\delta + \delta_1}{1 + \delta + \delta_1}}.   \]
It is possible  to choose $\delta$ and $\delta_1$ in such a way that
\[ a_1(\delta_1) \frac{\delta + \delta_1}{1 + 2 \delta_1} < \varepsilon. \]
\hfill $\Box$

\begin{theo} Let $s_j = \varphi(e_j)$ $( j =1, \ldots, r)$. Then
the height zeta function $Z_{\Sigma}(\varphi)$ is holomorphic for
${\rm Re}(s_j) > 1$.  There exists an analytic continuation
of $Z_{\Sigma}(\varphi)$ to the domain
${\rm Re}(s_j) > 1- \delta$ such that the only
singularities of $Z_{\Sigma}(\varphi)$ in this domain are poles
of order $ \leq 1$  along the hyperplanes $s_j = 1$ $(j =1, \ldots, r)$.
\label{extension.m}
\end{theo}

{\em Proof.} By the Poisson formula,
\[ Z_{\Sigma}(\varphi) =
\frac{1}{{\rm vol}(T^1({\bf A}_K)/T(K)) } \sum_{\chi \in
(T({\bf A}_K)/T(K))^* } \hat{H}_{\S}(\chi, -\varphi). \]
Since $H_{\Sigma}(x, -\varphi)$ is ${\bf K}_T$-invariant,
we can assume that in the above formula $\chi$ runs
over the elements of the group ${\cal P}$ consisting of
characters of $T({\bf A}_K)$ which are trivial on
${\bf K}_T \cdot T(K)$.

Let $J$ be a subset of $I = \{ 1, \ldots, r \}$.
Denote by ${\cal P}_J$ the
subset of ${\cal P}$ consisting of all characters $\chi \in {\cal P}$
such that the corresponding Hecke character $\chi_j$ is trivial
if and only if $j \in J$. Then
\[ Z_{\Sigma}(\varphi) = \sum_{J \subset I}
Z_{\Sigma,J}(\varphi) \]
where
\[ Z_{\Sigma,J}(\varphi) =
\frac{1}{{\rm vol}(T^1({\bf A}_K)/T(K)) } \sum_{\chi \in {\cal P}_J}
\hat{H}_{\S}(\chi, -\varphi). \]

Consider the logarithmic space
\[ N_{{\bf R}, \infty} =
\prod_{ v \in {\rm Val}_{\infty}(K)} T(K_v)/T({\cal O}_v) \]
containing the full sublattice $T({\cal O}_K)/W(T)$
of ${\cal O}_K$-integral points of $T(K)$ modulo torsion.
Let $\chi_{\infty}$ be the restriction of a character
$\chi \in {\cal P}$ to $N_{{\bf R}, \infty}$. Then $\chi_{\infty}(x) =
e^{2\pi i < x,y_{\chi} >}$ where $y_{\chi}$ is
an element of the dual logarithmic space
\[ M_{{\bf R}, \infty} =
\prod_{ v \in {\rm Val}_{\infty}(K)} {\rm Hom}\,( T(K_v)/T({\cal O}_v),
{\bf R}).  \]
Moreover, $y_{\chi}$ belongs to the dual lattice
$(T({\cal O}_K)/W(T))^*
\subset M_{{\bf R}, \infty}$.

Let $u_1, \ldots, u_{r'}$ be a basis of the lattice
$T({\cal O}_K)/W(T) \subset N_{{\bf R}, \infty}$, $f_1, \ldots, f_{r'}$
the dual basis of the dual lattice $(T({\cal O}_K)/W(T))^*
\subset M_{{\bf R}, \infty}$. We extend
\[ e^{2\pi i < x, f_1>} , \ldots, e^{2 \pi i<x, f_{r'}>} \]
 to some adelic characters
$\eta_1, \ldots, \eta_{r'} \in {\cal P}$. Using  \ref{subgroups} and
the basis $\eta_1, \ldots, \eta_{r'}$, we
can extend
\[ \chi_{\infty} = \prod_{k = 1}^{r'}
e^{2\pi i a_k < x,f_k >}, \;\; a_k \in {\bf Z} \]
to a character
\[ \tilde{\chi} = \prod_{k =1}^{r'} \eta_i^{a_k} \in {\cal P} \]
such that $\tilde{\chi}_{\infty} = \chi_{\infty}$ and
$\chi \cdot \tilde{\chi}^{-1}$ is a character of the finite group ${\bf
cl}(T)$.

We fix a character of $\chi_c$ of ${\bf cl}(T)$. Denote by
${\cal P}_{J,\chi_c}$ the
set of all characters $\chi \in {\cal P}_J$ such that
$\chi \cdot \tilde{\chi}^{-1} = \chi_c$.
Then a character $\chi \in {\cal P}_{J,\chi_c}$ is uniquely defined by
its archimedian component $\chi_{\infty}$.

By  \ref{dmethod},
\[ Q_{\S} (\chi,-\varphi) = \prod_{ v \not\in S_{\infty}}
\hat{H}_{\S,v}(\chi_v, \varphi)
\prod_{i =1}^r L_{K_j}^{-1}(s_j, \chi_j) \]
is absolutely convergent Euler product for ${\rm Re}(s_j) > 1 -
\delta > 1/2$.

By \ref{estim},
\[ \prod_{j \not\in J} L_{K_j}(s_j, \chi_j) < C(\varepsilon)
 \prod_{j \not\in J} {\rm Im}\, (s_j)^{\varepsilon}  \]
${\rm Re}(s_j) > 1 - \delta$.

We apply \ref{l-estimation} to the archimedian Fourier
transform
\[ \hat{H}_{\S,\infty}(\chi, \varphi) =
\prod_{ v \in S_{\infty}}
\hat{H}_{\S,v}(\chi_v, \varphi).  \]

Then, by \ref{lconver},
\[ \sum_{\chi \in {\cal P}_{J,\chi_c}}  \hat{H}(\chi, -\varphi)
\prod_{j \in J} \zeta_{K_j}(s_j)^{-1} \]
is absolutely convergent for ${\rm Re}(s_j) > 1 - \delta$.

Therefore, we have obtained that
\[ Z_{\S,J} (\varphi) \prod_{j \in J} \zeta_{K_j}(s_j)^{-1} \]
is a holomorphic function for ${\rm Re}(s_j) > 1 - \delta$
and for any $J \subset I$.

It remains to notice, that in the considered domain
 $\prod_{j \in J} \zeta_{K_j}(s_j)$
has only poles of order $1$ along hyperplanes $s_j = 1$.
\hfill $\Box$

\subsection{Rational points of bounded height}

Recall the standard tauberian statement:

\begin{theo}\cite{delange}
Let $X$ be a countable set, $F\, : \, X \ra {\bf R}_{>0}$ a real valued
function. Assume that
\[ Z_F(s) = \sum_{x \in X} F(x)^{-s} \]
is absolutely convergent for ${\rm Re}(s) > a> 0$ and has a representation
$$
Z_F(s)=(s-a)^{-r}g(s) + h(s)
$$
with $g(s)$ and $h(s)$ holomorphic for $Re(s)\ge a$, $g(a)\neq 0$,
$r\in N$.
Then for any $B>0$ there exists only a finite  number
$N(F,B)$ of elements $x \in X$ such that $F(x) \leq B$. Moreover,
\[ N(F,B) = \frac{g(a)}{a(b-1)!} B^a (\log B)^{r-1}(1+o(1)). \]
\label{tauberian}
\end{theo}

Let ${\cal L} = {\cal L}(\varphi_0)$ be a metrized invertible sheaf
over a smooth compactification ${\bf P}_{\Sigma}$ of an anisotropic
torus $T_K$  defined by a $G$-invariant fan $\Sigma$. We denote
by $Z_{\Sigma,{\cal L}}(s) = Z_{\Sigma}(s \varphi_0)$
the restriction of $Z_{\S}(\varphi)$ to the
line $s\lbrack {\cal L} \rbrack \subset {\rm Pic}({\bf P}_{\S})_{\bf R}$.
Let $a({\cal L})$ be the abscissa of convergence of $Z_{\Sigma,{\cal L}}(s)$
and
$b({\cal L})$ the order of the pole of $Z_{\Sigma,{\cal L}}(s)$ at
$s = a({\cal L})$.
By \ref{extension.m},
 \[ a({\cal L}) \leq \min_{j = 1}^r \frac{1}{\varphi(e_j)}. \]

By \ref{tauberian}, we obtain:

\begin{theo}
Assume that $\varphi(e_j) > 0$ for all $j =1, \ldots, r$; i.e.,
the class $\lbrack {\cal L} \rbrack$ is contained in the interior
of the cone of effective divisors $\Lambda_{\rm eff}(\Sigma)$.
Then there exists only finite number  $N({\bf P}_{\Sigma}, {\cal L},B)$ of
$K$-rational points $x \in T(K)$ having the
${\cal L}$-height $H_{\cal L}(x) \leq B$.
Moreover,
\[ N({\bf P}_{\Sigma}, {\cal L},B) =B^{a({\cal L})}\cdot (\log B)^{b({\cal
L})-1}(1+o(1)) \;\; ,
 B \ra \infty.  \]
\end{theo}

The following statement implies Batyrev-Manin conjectures
about the distribution of rational points of bounded ${\cal L}$-height
for smooth compactifications of anisotropic tori:

\begin{theo}
The number $a({\cal L})$ equals
\[ a({\cal L})={\rm inf} \, \{\lambda \mid \lambda \lbrack {\cal L}\rbrack  +
\lbrack {\cal K} \rbrack \in \Lambda_{\rm eff}(\Sigma)\};  \]
i.e.,
\[ a({\cal L}) = \min_{j = 1,r} \frac{1}{\varphi(e_j)}. \]
Moreover, $b({\cal L})$ equals the codimension of the minimal face
of  $\Lambda_{\rm eff}(\Sigma)$ containing
$a({\cal L})\lbrack {\cal L}\rbrack  +
\lbrack {\cal K} \rbrack$.
\end{theo}

{\em Proof. } By \cite{ono1}, we can choose the finite set $S$ such
that the natural homomorphism
\[ \pi_S \, : \, T(K) \rightarrow
\prod_{v \not\in S} T(K_v)/T({\cal O}_v)  =
\prod_{v \not\in S} N_v \]
is surjective.  Denote by $T({\cal O}_S)$ the kernel of
$\pi_S$ consisting of all $S$-units in $T(K)$.
The group $T({\cal O}_S)/W(T)$  has the natural embedding
in the finite-dimensional space
\[ N_{S,{\bf R}} = \prod_{v \in S} T(K_v)/T({\cal O}_v) \otimes {\bf R} \]
as a full sublattice.

Let $\Delta$ be the fundamental domain of  $T({\cal O}_S)/W(T)$ in
$N_{S,{\bf R}}$. For any $x \in T(K)$, denote by $\overline{x}_S$ the image
of $x$ in $N_{S,{\bf R}}$. Define $\phi(x)$ to be the
element of $T({\cal O}_S)$ such that $\overline{x}_S - \phi(x) \in \Delta$.
Thus, we have obtained the mapping
\[ \phi_S \, : \, T(K) \ra T({\cal O}_S). \]
Define the new height function $\tilde{H}_{\S}(x,\varphi)$ on $T(K)$ by
\[ \tilde{H}_{\S}(x, \varphi) = H_{\S}(\varphi, \phi_S(x))
\prod_{v \not\in S} H_{\S,v}(x_v, \varphi). \]
\noindent
Notice the following easy statement:

\begin{lem}
Choose a compact subset ${\bf K} \subset {\bf C}^r$ such that
${\rm Re}(s_j) > \delta$ $( j =1, \ldots, r)$ for $\varphi \in {\bf K}$.
Then there  exist positive constants $C_1$, $C_2$ such that
\[ 0 < C_1 < \frac{\tilde{H}_{\S}(x, \varphi)}{H_{\S}(x, \varphi)} <
C_2, \; \mbox{\rm for}\;  \varphi \in {\bf K}, \; x \in T(K). \]
\label{compare}
\end{lem}
Define $\tilde{Z}_{\S}(\varphi)$ by
\[ \tilde{Z}_{\S}(\varphi) = \sum_{x \in T(K)} \tilde{H}_{\S}(x, -\varphi). \]
Then $\tilde{Z}_{\S}(\varphi)$ splits into the
product
\[  \tilde{Z}_{\S}(\varphi) =
\prod_{v \not\in S} \left( \sum_{z \in N_v} H_{\S,v}(z,-\varphi) \right)
\cdot \left( \sum_{ u \in T({\cal O}) } H_{\S}(u, -\varphi) \right). \]
By \cite{drax1}, the Euler product
\[ \prod_{j =1}^r \zeta_{K_j}(s_j) \prod_{v \not\in S}
\left( \sum_{z \in N_v} H_{\S,v}(z,-\varphi) \right) \]
is a holomorphic function without zeros for ${\rm Re}(s_j) > 1/2$.

On the other hand,
\[ \sum_{ u \in T({\cal O}) } H_{\S}(u, -\varphi) \]
is an absolutely convergent series nonvanishing  for
${\rm Re}(s_j) > 0$.
Therefore, $\tilde{Z}_{\S}(\varphi)$ has a meromorphic extension
to the domain ${\rm Re}(s_j) > 1/2$ where it has poles of order $1$
along the hyperplanes $s_j = 1$.

By \ref{compare} and \ref{tauberian}, $\tilde{Z}_{\S}(\varphi)$ and
$Z_{\S}(\varphi)$ must have the same poles in the domain
${\rm Re}(s_j) > 1 - \delta$.  Therefore,
$Z_{\S}(\varphi)$ has poles of order $1$ along the hyperplanes $s_j = 1$.
By taking the restriction of $Z_{\S}(\varphi)$ to the line
$\varphi = s \varphi_0$, we obtain
the statement.
\hfill $\Box$

\subsection{The residue at $s_j = 1$}

Recall the  definition of the Tamagawa number of
Fano varieties \cite{peyre}.
This definition immediately extends to arbitrary
algebraic varieties $X$  with a metrized
canonical sheaf  ${\cal K}$.

Let $x_1, \ldots, x_d$ be local analytic coordinates on $X$. They
define a homeomorphism $f\,: \, U \ra {K_v}^d$
in $v$-adic topology between an
open subset $U \subset X$ and  $f(U) \subset {K_v}^d$.
Let $dx_1 \cdots dx_d$ be the Haar measure on $K_v^d$ normalized
by the condition
\[ \int_{{\cal O}_v^d} dx_1 \cdots dx_d = \frac{1}{(\sqrt{\delta_v})^d} \]
where $\delta_v$ is the absolute different of $K_v$.
Denote by $dx_1 \wedge \cdots \wedge dx_d$ the standard
differential form on $K_v^{d}$. Then
$g = f^*(dx_1 \wedge \cdots \wedge dx_d)$ is a local analytic section of
the metrized canonical sheaf ${\cal K}$.
We define the local measure on $U$ by
\[ \omega_{{\cal K},v} =f^*( \| g(f^{-1}(x)) \|_v dx_1 \cdots dx_d) . \]
The adelic Tamagawa measure  $\omega_{{\cal K},S}$ is defined  by
\[ \omega_{{\cal K},S} = \prod_{v \in {\rm Val}(K)} \lambda_v^{-1}
\omega_{{\cal K},v} \]
where  $\lambda_v = L_v(1, {\rm Pic}(X_{\overline{K}});\overline{K}/K))$
if $v \not\in S$,
$\lambda_v = 1$ if $v \in S$.

\begin{dfn}
{\rm \cite{peyre}
Let $\overline{X(K)}$ be the closure of $X(K) \subset X({\bf A}_K)$
in the
direct product topology. Then the  {\em Tamagawa number} of $X$ is defined by
\[ \tau_{\cal K}(X) =
\lim_{s \ra 1} (s-1)^r L_S(s, {\rm Pic}(X_{\overline{K}});\overline{K}/K))
\cdot \int_{\overline{X(K)}} \omega_{{\cal K},S}. \]}
\end{dfn}

\begin{prop}
Let ${\cal K} = {\cal L}(-\varphi_{\Sigma})$ be the metrized
canonical sheaf on a toric variety ${\bf P}_{\Sigma}$. Then
the restriction of the $v$-adic
measure $\omega_{{\cal K},v}$ to $T(K_v) \subset
{\bf P}_{\Sigma}(K_v)$ coincides with the measure
\[ H_{\S,v}(x, -\varphi_{\Sigma}) \omega_{\Omega,v}. \]
\label{restriction}
\end{prop}

{\em Proof.} The rational differential $d$-form $\Omega$ is a
rational section of ${\cal K}$. By definition of the $v$-adic metric on
${\cal L}(-\varphi_{\Sigma})$, $H_{\S,v}(x, -\varphi_{\Sigma})$
equals the norm $\| \Omega \|_v$ of the $T$-invariant
section $\Omega$. This implies the statement.
\hfill $\Box$

\begin{prop}
One has
\[ \int_{\overline{T(K)}} \omega_{{\cal K},S} =
\int_{\overline{{\bf P}_{\S}(K)}} \omega_{{\cal K},S}. \]
\end{prop}

{\em Proof.} Since $\overline{{\bf P}_{\S}(K)} \setminus \overline{T(K)}$
is a subset of ${\bf P}_T({\bf A}_K) \setminus T({\bf A}_K)$,  it is
sufficient to prove that
\[ \int_{T(K_v)} \omega_{{\cal K},v} =
\int_{{\bf P}_{\S}(K_v)} \omega_{{\cal K},v}. \]
Since the measure $\omega_{{\cal K},v}$ is $T({\cal O}_v)$-invariant and
the stabilizer in $T({\cal O}_v)$ of any point $x \in {\bf P}_{\S}(K_v)
\setminus T(K_v)$ is uncountable, the $\omega_{{\cal K},v}$-volume of
${\bf P}_{\S}(K_v) \setminus T(K_v)$ is zero.
\label{2int}
\hfill $\Box$.

\begin{theo}
Let $\Theta(\Sigma,K)$ be the
the residue of the zeta-function $Z_{\Sigma}(\varphi)$ at
$s_1 = \cdots = s_r = 1$. Then
\[ \Theta(\Sigma,K)  = \alpha({\bf P}_{\Sigma})\beta({\bf P}_{\Sigma})
\tau_{\cal K}({\bf P}_{\Sigma}).  \]
\end{theo}

{\em Proof.} By the Poisson formula,
\[ Z_{\Sigma}(\varphi) =
\frac{1}{{\rm vol}(T^1({\bf A}_K)/T(K)) } \sum_{\chi \in
(T({\bf A}_K)/T(K))^* } \hat{H}_{\S}(\chi, -\varphi). \]
By \ref{extension.m},
the residue of $Z_{\Sigma}(\varphi)$ at
$s_1 = \cdots = s_r = 1$ appears from $Z_{\S,I}(\varphi)$
containing only the terms $\hat{H}_{\S}(\chi, -\varphi)$
such that $\chi_1, \ldots , \chi_r$ are trivial characters of
${\bf G}_m({\bf A}_{K_j})/{\bf G}_m(K_j)$ $(j = 1, \ldots, l)$; i.e.,
$\chi$ is a character of the finite group $A(T)$.
We apply again the Poisson formula to the finite
sum
\[ Z_{\S,I}(\varphi)
= \frac{1}{{\rm vol}(T^1({\bf A}_K)/T(K)) }
\sum_{\chi \in (A(T))^* } \hat{H}_{\S}(\chi, -\varphi). \]
Using \ref{weak}, \ref{weak1}, \ref{obstr},  we have
\[ Z_{\S,I}(\varphi)
= \frac{\beta({\bf P}_{\Sigma})}{i(T_K)\cdot{\rm vol}(T^1({\bf A}_K)/T(K)) }
\int_{\overline{T(K)}} H_{\S}(x,-\varphi) \omega^1_{\Omega,S}. \]
Notice that $\omega^1_{\Omega,S} =  \omega_{\Omega,S}$ for
anisotropic tori.

Now we assume that $\varphi(e_1) = \ldots = \varphi(e_r) = s$.
Our purpose is to compute the constant
\[ \Theta(\Sigma,K) = \lim_{ s \ra 1} (s-1)^r Z_{\S,I}(s\varphi_{\Sigma}).
\]
By \ref{tamagawa1}, \ref{tamagawa2},
\[ \Theta(\Sigma,K)  =
\frac{\beta({\bf P}_{\Sigma})}{h(T_K)} L_S^{-1}(1, T; E/K)
\lim_{ s \ra 1} (s-1)^r
\int_{\overline{T(K)}} H_{\S}(x,-s\varphi_{\Sigma}) \omega_{\Omega,S}. \]
Notice that $\overline{T(K)}$ contains $T(K_v)$ for $v \not\in S$.
Denote by $\overline{T(K)}_S$ the image of $\overline{T(K)}$ in
$\prod_{v \in S} T(K_v)$.
By \ref{loc-integ}, we have
\[ \int_{\overline{T(K)}} H_{\S,v}(x,-s\varphi_{\Sigma}) \omega_{\Omega,S} =
\]
\[ = \prod_{v \not\in S} \int_{T(K_v)} H_{\S,v}(x,-s\varphi_{\Sigma}) d\mu_v
\cdot
  \int_{\overline{T(K)}_S}
\prod_{v \in S} H_{\S}(x,-s\varphi_{\Sigma}) \omega_{\Omega,v} =  \]
\[ =  L_S(s, T;E/K) \cdot \prod_{v \not\in S}
\left( \sum_{k =0}^d
  \frac{{\rm Tr}\, \Phi_v(i)}{q_v^{ks}} \right) \cdot
  \int_{\overline{T(K)}_S}
\prod_{v \in S} H_{\S,v}(x,-s\varphi_{\Sigma}) \omega_{\Omega,v}. \]
\[ L_S^{-1}(1, T; E/K)  \omega_{\Omega,S} =
\prod_{v \in {\rm Val}(K)} \omega_{\Omega,v}. \]
By \ref{p-function} and \ref{integral.1},
\[ L_S^{-1}(s, {\rm Pic}({\bf P}_{\Sigma,E}); E/K)
\prod_{v \not\in S} \left( \sum_{k =0}^d
   \frac{{\rm Tr}\, \Phi_v(i)}{q_v^{ks}} \right) \]
has no singularity at $s =1$.
Moreover, by \ref{restriction},
\[ \prod_{v \not\in S}
L_S^{-1}(1, {\rm Pic}({\bf P}_{\Sigma,E});E/K) \left( \sum_{k =0}^d
   \frac{{\rm Tr}\, \Phi_v(i)}{q_v^{k}} \right)
= \prod_{v \not\in S} \int_{T(K_v)} \lambda_v^{-1} \omega_{{\cal K},v}.  \]
Therefore
\[ \Theta(\Sigma,K) =
\frac{\beta({\bf P}_{\Sigma})}{h(T_K)}
\lim_{ s \ra 1} (s-1)^r
L_S(s, {\rm Pic}({\bf P}_{\Sigma,E}); E/K)
\int_{\overline{T(K)}} \omega_{{\cal K},S}. \]
It remains to apply \ref{2int}.
\hfill $\Box$

Using \ref{tauberian},  we obtain:

\begin{coro}
Let  $T$ be an anisotropic
torus and ${\bf P}_{\Sigma}$ its smooth compactification  $($notice that we do
not need to assume
that ${\bf P}_{\Sigma}$ is a Fano variety$)$. Let $r$ be
the rank of ${\rm Pic}({\bf P}_{\Sigma,K})$.
Then the number $N( {\bf P}_{\Sigma},{\cal K}^{-1}, B)$
of $K$-rational points  $x \in T(K)$ having the anticanonical
height $H_{{\cal K}^{-1}}(x) \leq B$ has the asymptotic
\[ N({\bf P}_{\Sigma},{\cal K}^{-1}, B) = \frac{\Theta(\Sigma,K)}{(r-1)!}
\cdot B (\log B)^{r-1}(1+o(1)), \hskip 0,3cm B\ra \infty.\]
\end{coro}

\end{document}